\documentclass{amsart}
\newtheorem{lemma}{Lemma}[section]
\newtheorem{propos}[lemma]{Proposition}

\newtheorem{theorem}[lemma]{Theorem}

\newtheorem{corol}[lemma]{Corollary}
\theoremstyle{definition}

\theoremstyle{remark}

\newcommand{\nquad}{\kern -10pt}

\newcommand{\CE}{\hbox{{$\mathcal E$}}}
\newcommand{\CS}{\hbox{{$\mathcal S$}}}

\newcommand{\Dsl}{{\nabla\!\!\!\!/}}
\newcommand{\bdel}{{\bar\del}}

\newcommand{\vardelta}{\mathfrak{\Delta}}

\newcommand{\C}{\mathbb{C}}
\newcommand{\R}{\mathbb{R}}

\newcommand{\Z}{\mathbb{Z}}

\newcommand{\del}{\partial}

\newcommand{\h}{{\scriptstyle\frac{1}{2}}}
\newcommand{\extd}{{\rm d}}
\newcommand{\isom}{{\cong}}
\newcommand{\eps}{{\epsilon}}
\newcommand{\tens}{\mathop{\otimes}}

\newcommand{\Ad}{{\rm Ad}}
\newcommand{\id}{{\rm id}}
\newcommand{\<}{\langle}
\renewcommand{\>}{\rangle}

\newcommand{\note}[1]{}
\newcommand{\vol}{\Upsilon}
\renewcommand{\o}{{}_{\scriptscriptstyle(1)}}
\renewcommand{\t}{{}_{\scriptscriptstyle(2)}}

\renewcommand{\th}{{}_{\scriptscriptstyle(3)}}
\newcommand{\fo}{{}_{\scriptscriptstyle(4)}}
\newcommand{\lvec}[1]{\overset{\leftarrow}{#1}}

\newcommand{\lbiprod}{{>\!\!\!\triangleleft\kern-.33em\cdot}}
\newcommand{\rbiprod}{{\cdot\kern-.33em\triangleright\!\!\!<}}

\newcommand{\eproof}{$\quad \diamond$\bigskip}

\newcommand{\eqn}[2]{\begin{equation}#2\label{#1}\end{equation}}

\begin{document}



\title[Riemannian geometry of the standard $q$-sphere]{\rm\large
NONCOMMUTATIVE RIEMANNIAN AND SPIN GEOMETRY OF THE STANDARD
$q$-SPHERE}

\author{S. Majid}%

\address{School of Mathematical Sciences\\
Queen Mary, University of London\\ 327 Mile End Rd,  London E1
4NS, UK}

\thanks{S.M. is a Royal Society University Research Fellow}%

\date{July 2003/Rev. 1}

\maketitle

\begin{abstract} We study the quantum sphere $\C_q[S^2]$ as a quantum
Riemannian manifold in the quantum frame bundle approach. We
exhibit its 2-dimensional cotangent bundle as a direct sum
$\Omega^{0,1}\oplus\Omega^{1,0}$ in a double complex. We find the
natural metric, volume form, Hodge * operator, Laplace and Maxwell
operators and projective module structure. We show that the 
$q$-monopole as spin connection
induces a natural Levi-Civita type connection and find its Ricci
curvature and $q$-Dirac operator $\Dsl$. We find the possibility
of an antisymmetric volume form quantum correction to the Ricci
curvature and Lichnerowicz-type formulae for $\Dsl^2$. We also
remark on the geometric $q$-Borel-Weil-Bott construction.
\end{abstract}

\section{Introduction}

The standard quantum sphere is nothing other than the invariant
subalgebra of the standard quantum group coordinate ring
$\C_q[SL_2]$ under a coaction of $\C[t,t^{-1}]$. In the
$*$-algebra setting it means $SU_2/U(1)$ in a coordinate form and
of course $q$-deformed. Other nonstandard quantum spheres were
constructed and classified in \cite{Pod:sph} while a unique
left-covariant 2-dimensional differential calculus on the standard
$\C_q[S^2]$ was found in \cite{Pod:dif,Pod:cla}. Meanwhile, the
$q$-monopole principal bundle with total space $\C_q[SL_2]$ and
base $\C_q[S^2]$ (i.e. the Hopf fibration) was constructed as an
example of the theory of quantum principal bundles in
\cite{BrzMa:gau} and has been somewhat studied since, so that many of
the ingredients of geometry for $\C_q[S^2]$ are already known.

In this paper we extend this geometry of $\C_q[S^2]$ to include
Riemannian structures and a geometrically natural Dirac operator,
using a systematic frame bundle approach to noncommutative
geometry in \cite{Ma:rie,Ma:rieq,Ma:ric}. In fact the point is not
just to obtain good geometrically justified proposals for these
structures on the algebra $\C_q[S^2]$ in isolation, but rather to
demonstrate that the frame bundle formulation, which is formulated
in principle at the level of any unital algebra, indeed includes
such an important example and gives reasonable answers for it.
This is important because without the straight-jacket of a general
theory that applies across diverse examples (including ones not
related to $q$-deformations) one could not have confidence that a
given definition was not ad-hoc, without which one could not
attach weight to physical or other predictions. We find that
$\C_q[S^2]$ indeed fits perfectly into this quantum frame bundle
approach to noncommutative geometry as a `quantum framed
Riemannian manifold'. Another motivation comes from the
operator-algebras and K-theory approach to noncommutative geometry
of Connes and others [Co]. How this can be reconciled with quantum
groups is an active and important area of research at the moment.
While we proceed only from the quantum groups side, by going up to
the point of the $q$-geometrically natural Dirac operator and
$q$-spin bundle, one can begin to compare with the `top down'
Connes approach where an axiomatically defined `Dirac' operator
implicitly defines the geometry. We find for example that our
$\Dsl$ indeed generates the exterior derivative by commutator, as
it should (equation (\ref{Dslcom})).

The principal features of the geometry on $\C_q[S^2]$ that we find
are as follows. The most important feature is that, unlike
previous examples based on quantum groups, the sphere is not
parallelizable. Hence there is no global `vielbein' and one must
work `upstairs' on the total space of the frame bundle for global
formulae. Unlike usual formulae in physics, in the noncommutative
case we do not consider coordinate charts or patching
transformations, instead we use only global constructions. Here
$\C_q[S^2]$ is the simplest example with this difficulty and hence
a good setting in which to demonstrate that the frame bundle
theory works. In fact we show, Theorem~2.1, that any quantum
homogeneous space induced by a Hopf algebra surjection is a framed
quantum manifold (we construct the soldering form). The
construction at the level of universal calculi was in
\cite{Ma:rie} but we extend this to general differential calculi
as needed for our example. To be self-contained, Section~2 starts
by recalling the required quantum homogeneous bundle
construction itself (before using it as frame bundle).

The rest of the paper computes what the general frame bundle
approach implies for the particular example of $\C_q[S^2]$. Here
for quantum frame bundle we take the quantum Hopf fibration or
`$q$-monopole principal bundle' from \cite{BrzMa:gau,BrzMa:rev},
in an appropriate form. The role of the fiber
 $SO(2)$ frame rotations is played by the commutative Hopf 
algebra $\C[t,t^{-1}]$ equipped with a noncommutative $q$-deformation 
of its usual
calculus. We then find the cotangent bundle $\Omega^1(\C_q[S^2])$
as an associated bundle to this. This means that, like all
associated bundles in this context (when there is a connection),
the cotangent bundle is necessarily projective, a point of view in
keeping with other approaches such as \cite{Con}. We later (in Section~5) 
exhibit its nontrivial projector explicitly. This is in spite of the
fact that we find (Theorem~3.1) that the cotangent bundle is
the sum of a charge -2 and charge 2 monopole, which means that is
is zero in the (noncommutative) K-theory of $\C_q[S^2]$ (this is actually in
keeping with the classical geometry). More importantly for us, Theorem~3.1 
 implies a natural direct
sum decomposition of the cotangent bundle into much simpler
``holomorphic'' and ``antiholomorphic'' parts
$\Omega^1(\C_q[S^2])=\Omega^{0,1}\oplus\Omega^{1,0}$ according to
the monopole charge, which we then use extensively in the sequel.

Section~4 covers the next `layer' of geometry in the form of the
exterior algebra, metric, Hodge
* operator, Laplace operator and Maxwell theory. In particular,
the natural $\C_q[SL_2]$-covariant metric $g$, Hodge $*$ and
volume 2-form (or symplectic structure) $\vol$ lifted to an
element $i(\vol)\in \Omega^1\bar\tens\Omega^1$ are naturally
related as $(\id\bar\tens *)(g)\propto i(\vol)$ (Proposition~4.3).

Section~5 then comes to the Riemannian geometry and contains our
main result (Theorem~5.1) that the $q$-monopole as spin connection
on the frame bundle indeed induces the correct generalized
Levi-Civita connection $\nabla$ on the cotangent bundle for its
natural metric. This is torsion free and `cotorsion free'. The
latter is a natural formulation of metric compatibility in the
skew form
\[ (\nabla\wedge\id-\id\wedge\nabla)g=0.\]
This condition was proposed in the axioms of \cite{Ma:rie,Ma:rieq} for
`quantum Riemannian manifolds' as a weakening of usual metric
compatibility suggested by noncommutative geometry. We see that
$\C_q[S^2]$ indeed bears this out. We also compute the Riemann
curvature of $\nabla$ as a 2-form-valued operator on 1-forms, and
from this the physically all-important Ricci curvature. There is
some freedom in the definition of this involved in choosing the
lifting map $i$ but we find that with the natural $i(\vol)$ modified by the
addition of a $q$-symmetric metric term, one has Ricci
proportional to the metric. Thus $\C_q[S^2]$ can be made into an
`Einstein space' (Proposition~5.2). Alternatively, we could not
modify $i(\vol)$, in which case we find
\[ {\rm
Ricci}={q^{-1}(1+q^4)\over 2}g+{[2]_q(1-q^4)\over 2}i(\vol)\]
showing a quantum correction involving the $q$-antisymmetric
volume form or symplectic structure. This effect could also be
formulated as a $q$-antisymmetric addition to the metric itself,
which would be in keeping with ideas from string theory, for
example. While the correct physical point of view and consequent
predictions of $q$-modifications to gravity would need an
understanding of the noncommutative stress energy tensor (which
will be attempted elsewhere), we see the possibility of a new
physical effect that vanishes as $q\to 1$.

Finally, for spin bundle we take $\CS=\CS_-\oplus \CS_+$ the
direct sum of the charge -1 and charge 1 $q$-monopole bundles.
This is then the correct `double cover' of the cotangent bundle in
terms of the corepresentation of $\C[t,t^{-1}]$. The same
$q$-monopole connection on the frame bundle as used for $\nabla$
now induces a covariant derivative $D$ on $\CS$. This combines
with a natural $\C_q[SL_2]$-covariant $\gamma$-matrix which we
provide, to give our gravitational Dirac operator $\Dsl$. It has
the correct $\Z_2$-graded form and we show (Proposition~5.5 and
equation (\ref{Dslsqac})) that its square is linked to the scalar
Laplacian. It is also relatively computable. For example,
\[ \begin{pmatrix} \pm q^{\h}a\\ b\end{pmatrix},\quad
\begin{pmatrix}\pm q^{\h}c\\ d\end{pmatrix}\]
are eigenspinors of mass $\pm q^\h$, where  $a,b,c,d$ are the
usual quantum group $\C_q[SL_2]$ generators, viewed now as spinor
components. We work algebraically and do not provide Hilbert space
or other analytic structures; this would need further study. We do, however,
show that unlike the  cotangent bundle, our spinor bundle
is trivial and we exhibit its trivialisation $\CS\isom \C_q[S^2]\oplus \C_q[S^2]$. 
Our $\Dsl$ in the trivialisation  appears
 to be more complicated than previous attempts at the
Dirac operator on the $q$-sphere such as \cite{PS:dir}, but comes
with the full geometrical picture above.

The appendix applies the $q$-monopole connection to formulate the
$q$-Borel-Weil-Bott construction as a byproduct of the
$q$-geometry in the paper. The generalisation of this to other
quantum groups and of the Riemannian and spin $q$-geometry to
other $q$-symmetric spaces are two directions for further work.
Notably for physics, a suitable $\C_q[S^4]$ and $q$-instanton are
known\cite{BCT:ins} but require the more general coalgebra bundle
theory for which nonuniversal differential calculi are not yet
formulated.

\subsection*{Acknoweldgements} I would like to thank Ruibin Zhang for
stimulating discussions on the Borel-Weil-Bott construction (see
appendix) during a visit to the Dept of Mathematics, University of
Sydney in December 2002.

\subsection*{Preliminaries}

We take $\C_q[SL_2]$ in the conventions of, for example, the text
\cite{Ma:book}. Namely it has a matrix $\begin{pmatrix}a&b\\\
c&d\end{pmatrix}$ of generators with $ba=qab$ etc. These are the
`lexicographical conventions' whereby $q$ is needed to put things
in lexicographical order. We will frequently use the
$q$-determinant relations $ad=1+q^{-1}bc$ and $da=1+qbc$. The Hopf
algebra structure has the usual matrix coproduct $\Delta$ and
counit $\eps$ on the generators and the antipode or `linearised
inverse' is $Sa=d$, $Sd=a$, $Sb=-qb$, $Sc=-q^{-1}c$. For the
axioms of Hopf algebra and basic notions such as actions and
coactions, we refer to \cite{Ma:book}. We use the Sweedler
notation\cite{Swe:hop} whereby $\Delta a=a\o \tens a\t$ and
$(\id\tens \Delta)\Delta a=a\o\tens a\t\tens a\th$, etc. We will
frequently need the right adjoint coaction $\Ad_R(a)=a\t \tens
(Sa\o)a\th$. We will write $A^+\subset A$ to denote the
augmentation ideal (the kernel of $\eps$). We take $q\in \C$
invertible. General constructions work in fact over any field but
our motivating point of view is over $\C$, which we retain for
convenience.

We recall that a differential calculus of an algebra $A$ means an
$A-A$-bimodule $\Omega^1$ and a map $\extd:A\to \Omega^1$ obeying
the Leibniz rule and such that $\Omega^1$ is spanned by 1-forms of
the form $a\extd b$. A calculus on a Hopf algebra is
left-covariant if the coproduct $\Delta:A\to A\tens A$ viewed as a
left coaction extends to a left coaction $\Delta_L$ on $\Omega^1$
such that $\extd$ is an intertwiner and $\Delta_L$ is a bimodule
map. In this case, having a Hopf-module, one knows that
$\Omega^1\isom A\tens \Lambda^1$ where $\Lambda^1$ are the
left-invariant 1-forms, and that $\Lambda^1=A^+/I$ where $I$ is
some right ideal contained in $A^+$ \cite{Wor:dif}.

 On the Hopf algebra $P=\C_q[SL_2]$ we take the 3-d calculus of
\cite{Wor:dif}. In our conventions this has a basis
\[  e^-=d\extd b-q b\extd
d,\quad  e^+=q^{-1} a\extd c-q^{-2}c\extd a,\quad  e^0=d\extd a-q
b\extd c\] left-invariant 1-forms, is spanned by these as a left
module (according to the above) while the right module relations
and exterior derivative are given in these terms by: \eqn{3dcom}{
e^\pm
\begin{pmatrix}a&b\\ c &
d\end{pmatrix}=\begin{pmatrix}qa&q^{-1}b\\ qc &
q^{-1}d\end{pmatrix}e^\pm,\quad  e^0
\begin{pmatrix}a&b\\ c &
d\end{pmatrix}=\begin{pmatrix}q^2 a&q^{-2}b\\ q^2c &
q^{-2}d\end{pmatrix} e^0\ }
\[ \extd a=a e^0+q b e^+,\quad \extd
b=a e^--q^{-2}b e^0,\quad \extd c=c  e^0+q d e^+,\quad \extd d=c
e^--q^{-2}d e^0.\] Our conventions for $e^\pm$ have been chosen
with hindsight to fit the frame bundle geometry, see
Theorem~\ref{Sman}. The corresponding ideal is \eqn{IP}{
I_P=\<a+q^2 d-(1+q^2),b^2,c^2,bc,(a-1)b, (d-1)c\>.}

Next, we let $A=\C[t,t^{-1}]$ be a Hopf algebra with $\Delta
t=t\tens t$ and $S t=t^{-1}$. This coacts on $\C_q[SL_2]$, making
it a comodule-algebra. Actually, a coaction here is the same thing
a $\Z$-grading and in our case the degrees are
\[\deg(a)=\deg(c)=1,\quad \deg(b)=\deg(d)=-1.\]
By definition the standard $q$-sphere $\C_q[S^2]$ is the degree
zero (i.e. invariant) subalgebra of $\C_q[SL_2]$. It is a
polynomial algebra $\C\<b_0,b_\pm\>$ with inherited relations
\[ b_\pm b_0=q^{\pm 2}b_0 b_\pm
\quad q^2b_- b_+=q^{-2}b_+ b_-+(1-q^{-2})b_0\] \eqn{Srel}{
b_0(1+qb_0)= b_+ b_-.} This last can also be written as
$b_0(1+q^{-1}b_0)=q^2b_-b_+$. Here $b_0=bc$, $b_+=cd$ and
$b_-=ab$. The first line of relations become as $q\to 1$ that the
algebra is commutative, while (\ref{Srel}) becomes the sphere
relation in term of $b_\pm$ complex and $b_0+1/2$. Moreover, the
coproduct of $\C_q[SL_2]$ restricts to $\C_q[S^2]$ as a left
coaction $\Delta_L:\C_q[S^2]\to \C_q[SL_2]\tens \C_q[S^2]$.
General (2-parameter) `quantum spheres' from the point of view of
left comodule algebras were obtained in \cite{Pod:sph}.

On this $q$-sphere we inherit a differential calculus from the one
above. It is not free over $\C_q[S^2]$ so we do not have a basis.
But it is spanned by \[\extd b_+=\extd(cd)=d^2 e^++c^2 e^-,\quad
\extd b_-=\extd(ab)=b^2 e^++a^2 e^-\] \eqn{db}{ \extd
b_0=\extd(bc)=q bd e^++q ac e^-} using the Leibniz rule and the
relations above. The inherited bimodule structure is far from
trivial and will be recovered below by our own means. It is
equivalent to formulae in \cite{Pod:dif}. The calculus inherits a
left coaction of $\C_q[SL_2]$ extending its coaction on
$\C_q[S^2]$.

Finally, the calculi in both cases extend to entire exterior
algebras. For $\C_q[SL_2]$ the natural extension compatible with
the super-Leibniz rule on higher forms and $\extd^2=0$ is:
\[ \extd  e^0=q^3 e^+\wedge e^-,\quad
\extd e^\pm=\mp q^{\pm 2}[2;q^{- 2}]e^\pm\wedge e^0,\quad
(e^\pm)^2=( e^0)^2=0\]
\[
q^2 e^+\wedge e^-+ e^-\wedge e^+=0,\quad  e^0\wedge e ^\pm+q^{\pm
4}e^\pm\wedge e^0=0\] where $[n;q]=(1-q^n)/(1-q)$ denotes a
$q$-integer. This means that there are the same dimensions as
classically, including a unique top form $ e^-\wedge e^+\wedge
e^0$. Again, these facts are well-known, but given here in our
required conventions. For $\C_q[S^2]$ the exterior calculus is not
so well-known and we obtain it below.

\section{Framings on nonuniversal quantum homogeneous spaces}

The general formulation of a quantum principal bundle with
nonuniversal calculi is as follows \cite{BrzMa:gau,BrzMa:rev}. As
`total space coordinate ring' we have an algebra $P$. For the
fiber a Hopf algebra $A$. We suppose that $P$ is a right
$A$-comodule algebra by a coaction $\Delta_R$ and define the fixed
subalgebra
\[ M=P^A=\{p\in P|\, \Delta_Rp=p\tens 1\}\]
for the `functions' on the base. For a bundle at the topological
level we require that  \eqn{bun0}{ 0\to P(\Omega^1M)P\to
\Omega^1P{\buildrel {\rm ver}\over \longrightarrow} P\tens A^+\to
0} is exact where $\Omega^1 P\subset P\tens P$ is the universal
calculus associated to any unital algebra (given by the kernel of
the product map). The map on the right is ${\rm ver}(p\tens
p')=p\Delta_R p'$, the generator of vertical vector fields. This
exactness is equivalent to the similar map $P\tens_M P\to P\tens
A$ being an isomorphism (a `Hopf-Galois' extension\cite{Sch:pri}).
When we have more general nonuniversal calculi then we require in
addition that \eqn{bun1}{ N_M=N_P\cap \Omega^1M}
\eqn{bun2}{\Delta_R N_P\subseteq N_P\tens A} \eqn{bun3}{{\rm
ver}(N_P)= P\tens I_A} where $\Omega^1(P)=(\Omega^1 P)/N_P$ etc
defines the calculus on $P$ as the quotient of the universal one
by a subbimodule. Here (\ref{bun1}) ensures that
\[ \Omega^1(M)={\rm span}\{m\extd_P n|\ n,m\in M\}\subseteq \Omega^1(P)\]
while (\ref{bun2}) ensures that $\Omega^1(P)$ is left covariant.
The coaction on $\Omega^1P$ here is the tensor product of the
coaction on each $P$. Finally, (\ref{bun3}) ensures that
\[ {\rm ver}:\Omega^1(P)\to P\tens \Lambda^1,\quad
\Lambda^1=A^+/I_A\] is well-defined (the calculus on $A$ is
assumed to be left-covariant as explained in the Preliminaries)
and exactness of \eqn{exact}{ 0\to P\Omega^1(M)P\to
\Omega^1(P){\buildrel {\rm ver}\over \longrightarrow} P\tens
\Lambda^1\to 0.} This is equivalent to the original formulation in
\cite{BrzMa:gau} based on such an exact sequence, as explained in
\cite{BrzMa:rev}. In effect, we put differential structures and
ensure that all relevant maps are compatible. Finally, in this
theory, a connection is defined\cite{BrzMa:gau} as an equivariant
splitting of $\Omega^1(P)$ providing a complement to the
`horizontal forms' $P\Omega^1(M)P$. If we assume that the calculus
on $A$ is bicovariant then a connection is equivalent to an
intertwiner
\[ \omega:\Lambda^1\to \Omega^1(P)\]
such that ${\rm ver}\circ\omega=1\tens\id$. Here $\Lambda^1$ has
the right adjoint coaction inherited from that on $A^+$.

For the purposes of this paper, the main example is a `quantum
homogeneous bundle' \cite{BrzMa:gau,BrzMa:rev} based on a
surjection $\pi:P\to A$ of Hopf algebras. This corresponds
geometrically to an inclusion of groups, and just as the subgroup
then acts by right multiplication, here $A$ coacts on $P$ by
\[ \Delta_R=(\id\tens \pi)\Delta:P\to P\tens A.\]
As above, we first construct the bundle with universal calculus
and `quantum homogeneous space base' $M=P^A$, and then impose
differential structures. To do the latter we assume that
$\Omega^1(P)$ is left-covariant and $\Omega^1(A)$ is bicovariant.
We take (\ref{bun1}) as a definition of $\Omega^1(M)$ while the
remaining conditions (\ref{bun2})-(\ref{bun3}) for a bundle with
these nonuniversal calculi reduce to \eqn{qhomog}{(\id\tens
\pi)\Ad_R(I_P)\subseteq I_P\tens A,\quad \pi(I_P)=I_A.} This
follows immediately using $N_P=\{p Sq\o\tens q\t|\ p\in P,\ q\in
I_P\}$ and computing $\Delta_R,{\rm ver}$ on such elements.

If one wants a connection, one can do this at the universal level
via a bicovariant splitting map $i:A\to P$. Thus,
\eqn{iomega}{\Delta_R \circ i=(i\tens \id)\Delta,\quad
(\pi\tens\id)\Delta \circ i=(\id\tens i)\Delta\quad \Rightarrow
\quad \omega(a)=Si(a)\o\extd i(a)\t} is a connection. One in fact
needs only a weaker $\Ad_R$-covariance condition \cite{BrzMa:gau}
but the stronger bicovariance implies this\cite{HajMa:pro} and is
the condition that is relevant below. In either case the map $i$
descends and defines a connection on the general bundle with
nonuniversal calculus if \eqn{ismooth}{ i(I_A)\subseteq I_P.} A
further refinement of these constructions for quantum principal
bundles can be found in \cite{BrzMa:rev}.

\note{\begin{eqnarray*} \Delta_R(p Sq\o\extd q\t)\nquad&& =p\o
Sq\o\t \tens q\t\o\tens \pi(p\t
(Sq\o\o)q\t\th)\\
&&=p\o Sq\t\o \tens q\t\t\tens \pi(p\t (Sq\o)q\th)\in N_P\tens A
\end{eqnarray*}
provided we assume that $q\t\tens \pi(Sq\o)q\th)\in I_P\tens A$
for all $q\in I_P$. This is the $\Ad_R$ condition stated. We used
coassociativity of the coproduct. Similarly
\[ {\rm ver}(p(Sq\o)\tens q\t)=p(Sq\o)q\t\o\tens \pi(q\t\t)=p\tens
\pi(q)\in P\tens I_A\] by the second condition stated. We see that
(\ref{bun3}) holds. We used cancellation via the axioms of the
antipode $S$.}

Up till now we have recalled the known quantum bundle construction
itself. We are now ready to give the nonuniversal version of the
frame bundle construction. An algebra $M$ is framed if it is the
base of a quantum principal bundle as above to which $\Omega^1(M)$
is an associated bundle. The frame quantum group fiber need not be
unique but its choice determines what kind of connections $\nabla$
on $\Omega^1(M)$ may be induced from connections on the frame
bundle. More details are in \cite{Ma:rie,Ma:rieq}. Apart from a
bundle over $M$ as above, we need an $A$-comodule $V$. Then
$\CE=(P\tens V)^A$ (the fixed submodule) plays the role of
sections of the associated bundle. Finally, we need a `soldering
form' $\theta:V\to P\Omega^1(M)$ such that the induced left
$M$-module map
\[ s_\theta:\CE\to \Omega^1(M),\quad p\tens v\mapsto p\theta(v)\]
is an isomorphism.

\begin{theorem}\label{homfra} Let $\pi:P\to A$ be a quantum homogeneous
bundle with general differential calculi as above. Then $M=P^A$ is
framed by the bundle and
\[ V=P^+\cap M/I_P\cap M,\quad \Delta_Rv=\tilde v\t\tens
S\pi(\tilde v\o)\]
\[ \theta(v)= S\tilde v\o\extd \tilde v\t\]
where $\tilde v$ is a representative of $v$ in $P^+\cap M$. Hence
every quantum homogeneous space of this type is a `quantum
manifold' in the framed sense.
\end{theorem}
\proof This construction for universal calculi is in
\cite[Prop.~4.3]{Ma:rie} so we have mainly to check that various
maps descend to the quotients needed for the nonuniversal calculi.
First observe that $v\in M$ means by definition $v\o\tens
\pi(v\t)=v\tens 1$. Moreover,  if $v\in M$ then $v\o\tens v\t\in
P\tens M$ because $v\o\tens v\t\o\tens \pi(v\t\t)=v\o\o\tens
v\o\t\tens \pi(v\t)=v\o\tens v\t\tens 1$. Similarly,
\[\Delta_Rv=v\o\tens \pi(Sv\o)=v\o\t\tens \pi(Sv\o\o)
\pi(v\t)=v\t\tens \pi(S v\o v\th)\] which is the projected adjoint
action. Hence if $v\in I_P\cap M$ we see from (\ref{qhomog}) and
from the above that $\Delta_Rv\in I_P\cap M\tens A$.  Hence
$\Delta_R$ descends to $V$. Incidentally, if $v\in P^+\cap M$ then
$\eps(v\t)\pi(Sv\o)=\pi(Sv)=S\pi(v\t)\eps(v\o)=1\eps(v)=0$ so
$\Delta_R$ is defined on $P^+\cap M$ in the first place (this is
the same as for the universal calculus case.) Meanwhile, if $v\in
I_P$ then $S\tilde v\o\tens \tilde v\t\in N_P$ and hence
$\theta(v)=0$ in $\Omega^1(P)$, so this is well-defined. Moreover,
if $\tilde v\in M$ is a representative of $v\in V$ then by the
above remark, $\theta(v)=S\tilde v\o \extd \tilde v\t\in
P\Omega^1(M)$ as required. That $\theta$ is equivariant follows
from this property proven for the universal calculi in
\cite{Ma:rie} to which we refer for the proof. Hence all maps are
defined as required and we have $s_\theta:(P\tens V)^A\to
\Omega^1(M)$. It remains to give its inverse, which we do by
quotienting the inverse in the universal calculus case, namely
\[ s_\theta^{-1}(m\extd n)=[m n\o\tens n\t-mn\tens 1],\quad \forall m,n\in M\]
where the expression in square brackets lies in $P\tens P^+\cap M$
(again using the observation above) and $[\ ]$ denotes the
equivalence class modulo $I_P\cap M$. That the result actually
lies in $(P\tens V)^A$ and gives the inverse of $s_\theta$ follows
in the same way as in the universal case in \cite{Ma:rie}. \eproof

Using $\omega$ such as from (\ref{iomega}), one may define the
covariant derivative
\[ D:\CE\to \Omega^1(M)\tens_M \CE,\quad D=(\id-\Pi_\omega)\extd\]
where we apply $\extd\tens\id$ to $\CE$ and
$\Pi_\omega=\cdot(\id\tens\omega)\circ {\rm ver}$ is the vertical
projection. When one takes the universal calculus this implies
that $\CE$ is a projective module and $D$ is the Grassmann
connection associated to the projector \cite{HajMa:pro}. In our
case we get other connections which we will study in the next
section. Also, using the framing, we obtain
\[ \nabla=(\id\tens s_\theta^{-1})\circ D\circ
s_\theta:\Omega^1(M)\to \Omega^1(M)\tens_M\Omega^1(M).\] Both $D$
and $\nabla$ behave as covariant derivatives (so
$\nabla(m\tau)=\extd m\tens_M\tau+m\nabla\tau$ for any function
$m\in M$ and 1-form $\tau$). Hence we need only give $\nabla$ on
exact forms.

\begin{propos} For the canonical connection induced by $i:A\to P$
on a quantum homogeneous space,
\[ \nabla(\extd m)=\extd(m\o S i\circ\pi(m\t)\o)\tens_M
i\circ\pi(m\t)\t Sm\th \extd m\fo.\]
\end{propos}

This is the nonuniversal version of a similar formula with the
universal differential calculus in \cite{Ma:rie}. The proof is
similar.

\section{Framing and holomorphic calculus on standard $q$-sphere $\C_q[S^2]$}

We start by recalling the known $q$-monopole bundle itself since
we will need it in full detail when we use it as frame bundle. We
fix the calculus on $P=\C_q[SL_2]$ to be the 3-d one as explained
in the Preliminaries. We take $A=\C[t,t^{-1}]$ and
\[ \pi(a)=t,\quad \pi(b)=\pi(c)=0,\quad \pi(d)=t^{-1}.\]
The right coaction $\Delta_R=(\id\tens\pi)\Delta$ works out as
\[ \Delta_R\begin{pmatrix}a&b\\ c&d\end{pmatrix}=\begin{pmatrix}a&b\\
c&d\end{pmatrix}\tens \begin{pmatrix}t& 0\\
0&t^{-1}\end{pmatrix}\] so that $\Delta_R a=a\tens t$ etc.
corresponding to $\deg(a)=1$, etc. Then $M=P^A=\C_q[S^2]$ as
explained in the Preliminaries. It is known that we have a quantum
homogeneous bundle with universal calculi. As in \cite{BrzMa:gau}
we then take
\[ I_A=\pi(I_P)=\<t+q^2 t^{-1}-(1+q^2)\>=\<(t-1)(t-q^2)\>\]
(we factored $t^{-1}$ out of the generator obtained from
projecting those of $I_P$). Now, this $I_A$ defines the
1-dimensional calculus on $\C[t,t^{-1}]$ with basis $\extd
t=t\tens [t-1]$ (where $[\ ]$ denotes modulo $I_A$) and relations
\[ \extd t. t=t^2\tens [t-1]t=t^2\tens [(t-1)(t-q^2)]+t^2\tens
q^2[t-1]=q^2t^2\tens[t-1]=q^2 t\extd t\] This is a
$q$-differential calculus whereby
\[ \extd (t^m)=[m;q^2]t^{m-1}\extd t\]
so that the relevant partial derivative is the usual
$q$-derivative.

We also verify that $I_P$ obeys the $\Ad_R$ condition in
(\ref{qhomog}). Indeed, the element $a+q^2 d$ (the $q$-trace) is
$\Ad_R$-invariant.  Meanwhile
\[ (\id\tens\pi)\Ad_R(b^2)=b^2\tens t^{-4},\quad
(\id\tens\pi)\Ad_R(c^2)=c^2\tens
t^{4}\] and so forth. Hence we have the quantum sphere as a
quantum homogeneous space where the calculus on it is obtained by
restriction of that on $\C_q[SL_2]$ as required in (\ref{bun1}).

Finally, \[i(t^n)=a^n,\quad i(t^{-n})=d^n,\quad \forall n\ge 0\]
defines a natural connection in the bundle via (\ref{iomega}) as
follows. We have
\[ i(I_A)={\rm
span}\{a^m(a-1)(a-q^2),\, a+q^2d-(1+q^2),\,
d^m(d-1)(d-q^2)\}\subseteq I_P.\] The middle term is already in
$I_P$. Hence also, multiplying it by $(a-q^2)$ we have
\[ (a-1)(a-q^2)+q^2(d-1)(a-q^2)\in I_P.\]
The second term is $q^2(da-1-q^2d-a+q^2+1)$ which lies in $I_P$
since $da-1=qbc\in I_P$. Hence $(a-1)(a-q^2)\in I_P$. Similarly
for $(d-1)(d-q^2)$. Hence the canonical connection defined by this
$i$ descends to the chosen nonuniversal calculi. The resulting
$q$-monopole connection is
 \eqn{omegamon}{ \omega(t^n)=
[n;q^2]e^0} for all integers $n$. This is easily proven by
induction as follows. Thus, when $n\ge 0$,
\begin{eqnarray*} \omega(t^n)\nquad &&=(Sa^n\o\extd a^n\t)
=S(a\o a'\o a''\o \cdots)\extd (a\t
a'\t a''\t\cdots)\\
&&=S(a'\o a''\t\cdots)Sa\o\left((\extd a\t)a'\t
a''\t\cdots+ a\t\extd(a'\t a''\t\cdots)\right)\\
&&=\omega(t^{n-1})+S(a'\o a''\t\cdots) \omega(t) a'\t a''\t \cdots\\
&&=\omega(t^{n-1})+S(a'\o a''\t\cdots) e^0 a'\t a''\t
\cdots=\omega(t^{n-1})+q^{2(n-1)} e^0
\end{eqnarray*}
where $a^n=a a' a''\cdots $ is the product of $n$ copies of the
generator $a\in \C_q[SL_2]$ (the primes are to keep the instances
apart). We used the antimultiplicativity of the antipode $S$ and
the Leibniz rules. We then used that $\omega(t)=(Sa\o)\extd
a\t=d\extd a-q b\extd c=e^0$ from the definition of the 3-d
calculus in the preliminaries. Finally, we used that $a'\t
a''\t\cdots$ has degree $n-1$ and hence its commutation relations
with $e^0$ give a factor $q^{2(n-1}$ after which we cancel using
the antipode axioms and $\eps(a)=1$. When we do the same with
$n=-n'$, $n'\ge 0$ giving $\omega(t^n)=Sd^{n'}\o\extd d^{n'}\t$
and a factor $q^{-2(n'-1)}$ at the corresponding point. We also
have $\omega(t^{-1})=-q^{-2}e^0$. The two halves of the
computation combine to the uniform answer (\ref{omegamon}). The
curvature of the $q$-monopole connection is \eqn{Fmon}{
F_\omega(t^n)=\extd\omega(t^n)+\omega(t^n)\wedge\omega(t^n)=[n,q^2]\extd
e^0=q^3[n;q^2]e^+\wedge e^-.} These constructions so far are not
essentially new. They are a version of the $q$-monopole
construction in \cite{BrzMa:gau,BrzMa:rev} . The choices and
conventions are slightly closer to those in \cite{HajMa:pro}
where, however, only universal calculi were considered.

Next, we compute $V$ in Theorem~\ref{homfra}. Clearly $P^+\cap
M=M^+=\ker\eps_M$. In our case $M=\C_q[S^2]$ is generated by $1,
b_\pm, b_0$ so that $M^+=\<b_0,b_\pm\>$ as an ideal. Meanwhile,
because $(a-1)b, (d-1)c, a+q^2d-(1+q^2)$ are not of homogeneous
degree, the ideal which each one generates has no intersection
with $M$. We assume here that $\C_q[SL_2]$ has no zero-divisors.
We therefore focus on $b^2,c^2,bc$. The elements of degree zero in
$\<b^2\>$ include $b^2\{a^2,ac,c^2\}$. Hence, $b_-^2,b_-b_0,b_0^2$
lie in $I_P\cap M$. Similarly from $\<c^2\>$ we have
$b_+^2,b_0b_+$ also in this ideal. The element $bc=b_0$ is already
in the ideal. From these considerations we arrive at
\[ V=\<b_\pm\>/\<b_\pm^2,b_0\>.\]
Hence $V$ is 2-dimensional with representatives $b_\pm$. We then
compute the coaction $\Delta_R$ on $V$ from Theorem~\ref{homfra}
as \eqn{coactV}{  \Delta_R b_+=cd\tens S\pi(d^2)=b_+\tens
t^2,\quad \Delta_R b_-=ab\tens S\pi(a^2)=b_-\tens t^{-2}.} Hence
$V=\C\oplus \C$ and the associated bundle \[ \CE=\CE_{-2}\oplus
\CE_{+2}=\C_q[SL_2]_{2}\oplus \C_q[SL_{2}]_{-2}\] is the direct
sum of the $q$-monopole bundles of charge -2 and charge 2. We
identify their sections with  the $\pm2$ degree components in
$\C_q[SL_2]$. Thus Theorem~\ref{homfra} yields:

\begin{theorem}\label{Sman} $\C_q[S^2]=\C_q[SL_2]_0$ is a framed
quantum manifold with cotangent bundle
 \[ \Omega^1(\C_q[S^2])\isom \CE_{-2}\oplus
\CE_{2}\] isomorphic to the charge 2 and charge -2 monopole
bundles. This isomorphism is given by the soldering form
\[ \theta(b_-)=d^2\extd b_-+q^2 b^2\extd
b_+-[2;q^2]bd\extd b_0= e^-\]
\[ \theta(b_+)=a^2\extd b_++q^{-2}c^2\extd
b_--[2;q^{-2}]ac\extd b_0=e^+\] and makes $\Omega^1(\C_q[S^2])$
projective. \end{theorem} \proof The only remaining part is to
compute $\theta(b_\pm)$. We first find the coaction on $\C_q[S^2]$
inherited from the coproduct of $\C_q[SL_2]$ as
\[
\Delta_L(b_-)=\Delta(ab)=ab\tens(1+[2]_qb_0)+a^2\tens b_-+b^2\tens
b_+\]
\[ \Delta_L(b_+)=\Delta(cd)=cd\tens (1+[2]_q b_0)+c^2\tens
b_-+d^2\tens b_+\] \eqn{coactS}{ \Delta_L(b_0)=\Delta(bc)=1\tens
b_0+bc\tens(1+[2]_qb_0)+qac\tens b_-+qbd\tens b_+} where
$[2]_q=q+q^{-1}$. These coproducts were already used in computing
$\Delta_R b_\pm$ above; this time we apply $S$ to the first factor
and compute
\[ \theta(b_+)=Sb_+\o\extd b_+\t=-q^{-1}ac\extd b_3+a^2\extd
b_++q^{-2}c^2\extd b_-+q(-q^{-1})ac\extd b_3 q.\] Similarly for
$\theta(b_-)$. This gives the middle expressions. We then insert
(\ref{db}) and find $e^\pm$ for the values of the map $\theta:V\to
\Omega^1(\C_q[SL_2])$. By Theorem~\ref{homfra} this is
well-defined on $V$ and actually has its values in
$\C_q[SL_2]\Omega^1(\C_q[S^2])$. Also according to
Theorem~\ref{homfra}, one must multiply $\theta(b_-)$ by an
element of degree 2, and $\theta(b_+)$ by an element of degree -2
to get 1-forms on $\C_q[S^2]$ and every 1-form is obtained in this
this way. $\CE_{\pm 2}$ are both projective as shown in
\cite{HajMa:pro}, as given via the Cuntz-Quillen theorem and the
$q$-monopole connection with the universal calculus. \eproof

These two 1-forms $\theta(b_\pm)=e^\pm$ play the role of
`vielbein' but do not themselves live on the base. Not also that
as regards the bimodule structure, the elements of $\C_q[S^2]$
commute with the $\theta(v)$ as we see from (\ref{3dcom}).

\begin{corol}\label{dSrel} The three 1-forms $\extd b_\pm$ and $\extd b_0$
enjoy the relation
\[ q^2 b_- \extd b_+ +  b_+ \extd b_- - (1+[2]_q b_0)\extd b_0=0\]
where $[2]_q=q+q^{-1}$ denotes a symmetrized $q$-integer.
\end{corol} \proof Using (\ref{coactS}) we have
$\theta(b_0)=Sb_0\o\extd b_0\t=b_0\extd (1+qb_0)+ (-q)q b_-\extd
b_++(-q^{-1})qb_+\extd b_-+(1+q^{-1}b_0)\extd b_0=0$ since $b_0$
represents zero in $V$. This identity can also be obtained from
requiring $\extd b_\pm$ and $\extd b_0$ to be recovered from
(\ref{db}) composed with Theorem~\ref{Sman} and extensive use of
the commutation relations. \eproof

This identity corresponds in the classical case to the
differential of (\ref{Srel}). However, for $q\ne 1$ this is not so
immediate because the bimodule relations for $\Omega^1(\C_q[S^2])$
are complicated to find explicitly. On the other hand,
Theorem~\ref{Sman} implies a direct sum structure to the cotangent
bundle with each piece more reasonable to work with.

\begin{corol}\label{dhol} $\Omega^1(\C_q[S^2])=\Omega^{0,1}\oplus
\Omega^{1,0}$ where $\Omega^{1,0},\Omega^{0,1}$ are each first
order left-covariant differential calculi over $\C_q[S^2]$ with
differentials $\del,\bdel$ obeying $\extd=\del+\bdel$ and
\[ \del b_-\begin{cases}b_+\\ b_-\\ b_0\end{cases}
\nquad=\begin{cases}q^{-2}b_+\del b_-\\ q^{2} b_-\del b_-\\
b_0\del b_-,\end{cases}\nquad\ \del b_+\begin{cases}b_+\\ b_-\\
b_0\end{cases}\nquad=\begin{cases}q^{2}b_+\del b_+\\
q^2b_-\del b_++(q^2-q^{-2})b_+\del b_-\\
q^4 b_0\del b_+\end{cases}\]
\[ \bdel b_-\begin{cases}b_+\\ b_-\\ b_0\end{cases}
\nquad=\begin{cases}q^{-2}b_+\bdel b_-+(q^{-2}-q^2)b_-\bdel b_+\\
q^{-2} b_-\bdel b_-\\
q^{-4}b_0\bdel b_-,\end{cases}\nquad \bdel b_+\begin{cases}b_+\\ b_-\\
b_0\end{cases}
\nquad=\begin{cases}q^{-2}b_+\bdel b_+\\ q^{2} b_-\bdel b_+\\
b_0\bdel b_+.\end{cases}\] One has the relations
\[ \del b_0=q^2
b_-\del b_+-q^{-2}b_+\del b_-,\quad \bdel b_0=b_+\bdel
b_--q^4b_-\bdel b_+\]
\[
b_0b_-\del b_+=q^{-3}(1+q^{-1}b_0)b_+\del b_-,\quad b_0b_+\bdel
b_-=q^3(1+qb_0)b_-\bdel b_+.\]
\end{corol}
\proof From Theorem~\ref{Sman} know that $\Omega^1$ is a direct
sum as a left module over $\C_q[S^2]$ spanned respectively by
\[
\{b^2,db,d^2\}e^+=\{\del b_-,\del b_0,\del b_+\},\quad
\{a^2,ca,c^2\}e^-=\{\bdel b_-,\bdel b_0,\bdel b_+\} \] where we
use the expressions (\ref{db}) to identify the two components of
$\extd$. Next we observe that $e^\pm$ commute with elements of
$\C_q[S^2]$ so that the commutation relations of functions with
$\del$ and $\bdel$ are easily determined from the relations of
$\C_q[SL_2]$. We find the ones stated and
\[ \del b_0\begin{cases}b_+\\ b_-\\
b_0\end{cases}\nquad=\begin{cases}b_+\del b_0\\
q^4b_-\del b_0+q(1-q^2)\del b_-\\
q^{2}b_0\del b_0,\end{cases}\quad \bdel b_0\begin{cases}b_+\\
b_-\\
b_0\end{cases}\nquad=\begin{cases}q^{-4}b_+\bdel b_0
+q^{-1}(1-q^{-2})\bdel b_+\\  b_-\bdel b_0\\
q^{-2}b_0\bdel b_0.\end{cases} \] These close so that each
$\del,\bdel$ generate a bimodule differential calculus. Their
Leibniz rules follow from that for $\extd$ and the direct sum
decomposition. Next, we observe the relations
\[  b_+\del b_-=qb_0\del b_0,\quad b_-\del
b_+=q^{-2}(1+q^{-1}b_0)\del b_0\] \[ b_-\bdel b_+=q^{-3}b_0\bdel
b_0, \quad b_+\bdel b_-=(1+qb_0)\bdel b_0\] following likewise
from the relations of $\C_q[SL_2]$ acting on $e^\pm$. We use them
as a definition of $\bdel b_0,\del b_0$ and a relation among the
$\bdel b_\pm$ (respectively, $\del b_\pm$) as stated. The latter
also imply the relation in (\ref{dSrel}) and follow from
differentiating (\ref{Srel}), which is their geometrical content
(related to the rank 1 projective module structure of each
bundle). We note that one also has other relations, such as
\[ b_+\del b_0=q^2b_0\del b_+,\quad b_-\del
b_0=q^{-1}(1+q^{-1}b_0)\del b_-\] \[ b_-\bdel b_0=q^{-2}b_0\bdel
b_-,\quad b_+\bdel b_0=q(1+qb_0)\bdel b_+\]   which are not
independent of the one already found. For example the  $\bdel$
relations here can be written as
\[b_-^2\bdel b_+=q^{-7} b_0^2\bdel b_-,\quad b_+^2\bdel
b_-=q(1+q^3b_0)(1+qb_0)\bdel b_+\] which can be deduced from the
one stated if one assumes that the left action of $b_0$ and
$1+q^{-1}b_0$ respectively can be cancelled. Finally, each of the
spaces $\Omega^{1,0}$ and $\Omega^{0,1}$ are stable under the left
coaction of $\C_q[SL_2]$. This is because `upstairs' on
$\Omega^1(SL_2)$ the coaction on an element $fe^\pm$ is just
$f\o\tens f\t e^\pm$ and the coproduct defines a left coaction in
each degree (because left-comultiplication commutes with the
right-comultiplication  used in defining the grading according to
$(\id\tens\pi)\Delta$). This coaction is intertwined by $\del$ and
$\bdel$ with the left coaction (\ref{coactS}) on $\C_q[S^2]$ since
this is true for $\extd$. \eproof

Clearly these `holomorphic' and `antiholomorphic' cotangent
bundles are relatively simple to work with. On the other hand, we
can use Theorem~\ref{Sman} to compute $\del,\bdel$ in terms of
$\extd$. Using the relations in $\C_q[SL_2]$ and
Corollary~\ref{dSrel} we find
\[ \del
b_-=qb_-\extd b_0-q^{-1}b_0\extd b_-,\quad \del b_+=(1+qb_0)\extd
b_+-q^{-1}b_+\extd b_0\]
\[\del
b_0=q^2b_-\extd b_+-q^{-1}b_0\extd b_0,\quad \bdel b_0=b_+\extd
b_- - qb_0\extd b_0\] \eqn{deld}{\bdel b_-=(1+q^{-1}b_0)\extd
b_--qb_-\extd b_0,\quad \bdel b_+=q^{-1} b_+\extd b_0-q b_0\extd
b_+.  } Therefore, as an application, we can recover the bimodule
relations in $\Omega^1(\C_q[S^2])$ from the much simpler ones for
the two parts.

\begin{propos}\label{dbbimod} Let $\mu=q^2-q^{-2}$. The bimodule
relations for the 2-dimensional calculus on $\C_q[S^2]$ are
\[ \extd b_0\begin{cases}b_0\\ b_\pm\end{cases}
\nquad=\begin{cases}(q^2+q\mu b_0) b_0\extd b_0-\mu b_0 b_+\extd b_-\\
q^{\mp2}(1\mp q^{\pm1}\mu b_0)b_\pm\extd b_0-(1-q^{\pm2}\mp q^{\pm
1}\mu b_0)b_0\extd b_\pm
\end{cases}\]
\[ \extd b_\pm\begin{cases}b_0\\ b_\mp\\
b_\pm\end{cases}\nquad =\begin{cases}(q^{\pm 4}\pm q^{\pm 1}\mu
b_0)b_0\extd b_+\mp q^{\mp 1}\mu b_\pm b_0\extd b_0\\ q^{\pm
2}(1\pm q^{\pm 1}\mu b_0)b_\mp\extd b_\pm \pm  q^{\pm 1}\mu
q^{-1}b_0^2\extd b_0\\ q^{\pm2}(1\pm q^{\pm 1}\mu b_0)b_\pm\extd
b_\pm\mp q^{\mp 1}\mu b_\pm^2\extd b_0.
\end{cases}\]
\end{propos}
\proof These are all computed along the following lines: \[\extd
b_0.b_0=q^{-2}b_0\bdel b_0+q^2b_0\del b_0=q^2b_0\extd
b_0+(q^{-2}-q^2)b_0\bdel b_0\] using $\extd=\del+\bdel$ and the
commutation relations from Corollary~\ref{dhol}. We then express
$\bdel b_0$ in terms of $\extd$ to obtain the result. Similarly
for all the other commutation relations.  \eproof

As a cross-check, one may now verify that Corollary~\ref{dSrel}
corresponds to the differentials of the the $q$-sphere relations.
Put another way, using Corollary~\ref{dSrel}, the differentials of
the four relations (\ref{Srel}) of the $q$-sphere reduce to
\[ q^{-1}\extd b_+. b_--\extd b_0. b_0= -q^{-2}b_0 \extd b_0+q
b_-\extd b_+\]
\[ q\extd b_-. b_+-q^{-2}\extd b_0. b_0= -b_0 \extd b_0+q^{-1}
b_+\extd b_-\] \[\extd b_\pm.b_0-q^{\pm 2}\extd b_0.b_\pm=q^{\pm
2}b_0\extd b_\pm-b_\pm \extd b_0\] which all hold using
Proposition~\ref{dbbimod}. In fact Podles in \cite{Pod:cla} has
shown that there is a unique left-covariant calculus on
$\C_q[S^2]$ of the correct classical dimension, hence
Proposition~\ref{dbbimod} is necessarily isomorphic to this, but
derived differently.

\section{Exterior algebra, Hodge-$*$ and Maxwell theory on the $q$-sphere}

In the last section we have expressed the cotangent bundle of
$\C_q[S^2]$ as associated to a frame bundle by a $1\oplus
1$-dimensional representation of the `frame group' $\C[t,t^{-1}]$
equipped with a (bicovariant) $q$-differential structure. We
deduced that $\Omega^1(\C_q[S^2])$ is the sum of two 1-dimensional
left-covariant calculi $\Omega^{0,1}$ and $\Omega^{1,0}$. We now
extend this to the entire exterior algebra.

First of all, working `upstairs' in $\Omega(\C_q[SL_2])$ we define
$\Omega(\C_q[S^2])$ as the differential algebra obtained by
restriction. It is generated by $\del b_\pm$ and $\bdel b_\pm$
which means generated by $e^\pm$ and certain elements of
$\C_q[SL_2]$. Because of the commutation relations with $e^\pm$,
these elements can all be collected to the left. Because of the
relations between the $e^\pm$, there is only a functional multiple
of $e^+\wedge e^-$ in degree 2 and nothing in higher degree, so
\[ \Omega^2(\C_q[S^2])=\Omega^{1,1},\quad \Omega^{2,0}=0=\Omega^{0,2}\]
where the numbers refer to the degrees in $\del,\bdel$ and
$\Omega^{1,1}$ is 1-dimensional over the algebra. We also extend
$\del$ and $\bdel$ by $\del^2=\bdel^2=0$ so that each generates an
exterior algebra (with top degree 1 by the above arguments
`upstairs' or from the relations in Corollary~\ref{dhol}. Finally,
we define
\[
\del=\extd|_{\Omega^{0,1}},\quad \bdel=\extd|_{\Omega^{1,0}}\]
giving a
double complex \begin{eqnarray*} 0& & \ 0\\
{\scriptstyle \bdel} \uparrow& &{\scriptstyle \bdel} \uparrow \\
 \Omega^{0,1}&{\buildrel
\del\over\longrightarrow}&\Omega^{1,1}\quad {\buildrel
\del\over\longrightarrow}\quad 0\\
{\scriptstyle \bdel} \uparrow & &{\scriptstyle \bdel} \uparrow\\
\Omega^{0,0}&{\buildrel
\del\over\longrightarrow}&\Omega^{1,0}\quad {\buildrel
\del\over\longrightarrow}\quad 0
\end{eqnarray*}
Here the graded-derivation property of $\extd$ implies that \[
\bdel\del+\del\bdel=0.\]

Moreover, forms which are left-invariant under the $\C_q[SL_2]$
coaction are precisely the ones generated by $e^\pm$ alone. Thus,
up to scale, there is a unique left-invariant top form
\[ \vol=e^+\wedge e^-.\]
This is a basis of $\Omega^2(\C_q[S^2])$ over $\C_q[S^2]$. We let
$\mu=q^2-q^{-2}$.

\begin{propos}\label{delbardel} The relations between the $\Omega^{0,1}$
and $\Omega^{1,0}$ calculi are
\[\del b_+\wedge\bdel
b_-+q^6\bdel b_-\wedge\del b_+=q^4\mu (b_0^2-1)\vol \]
\[ \del
b_-\wedge \bdel b_+=-q^2\bdel b_+\wedge \del b_-=q^2b_0^2\vol\]
\[ \del b_-\wedge \bdel b_-=-q^6\bdel b_-\wedge\del
b_-=q^5b_-^2\vol \]
\[ \del b_+\wedge\bdel b_+=-q^6\bdel b_+\wedge\del b_+=q^5
b_+^2\vol \]
\end{propos}
\proof We compute all expressions in terms of $e^\pm$ using the
definitions from the proof of Corollary~\ref{dhol}. For example
\[ \bdel b_-\wedge \del b_+=a^2 e^-\wedge d^2
e^+=q^{-2}a^2d^2e^-\wedge e^+=-(1+q^{-3}b_0)(1+q^{-1}b_0)\vol\]
using the relations between functions and $e^\pm$ in
$\Omega^1(\C_q[SL_2])$ and the relations in the quantum group.
Computing all expressions in this way and comparing gives the
relations stated. One may similarly compute
\[
\del b_0\wedge\bdel b_0=q^4(1+qb_0)b_0\vol,\quad \bdel
b_0\wedge\del b_0=-(1+q^{-1}b_0)b_0\vol\]
\[
\del b_0\wedge\bdel b_-=q^4(1+qb_0)b_-\vol,\quad \bdel b_-\wedge
\del b_0=-(1+q^{-3}b_0)b_-\vol\]
\[ \del b_+\wedge\bdel b_0=q^4b_+(1+q b_0)\vol,\quad \bdel
b_0\wedge\del b_+=-b_+(1+q^{-3}b_0)\vol\]\[ \del b_0\wedge\bdel
b_+=-q^4\bdel b_+\wedge \del b_0=q^3b_+b_0\vol,\quad  \del
b_-\wedge\bdel b_0=-q^4\bdel b_0\wedge\del b_-=q^5 b_- b_0\vol\]
which will be useful later on, giving the further relations
\[ q^{-4}\del b_0\wedge \bdel
b_0+\bdel b_0\wedge\del b_0=(q-q^{-1})b_0^2\vol\] \[ \del
b_0\wedge\bdel b_-+q^8\bdel b_-\wedge\del b_0=-q^6\mu
b_-\vol,\quad \del b_+\wedge\bdel b_0+q^8\bdel b_0\wedge\del
b_+=-q^6\mu b_+\vol.\hfil \mbox{\eproof}\]  \goodbreak

 Note that first two lines taken together exhibit $\vol$ as an element of
$\Omega^{1,1}(\C_q[S^2])$, which is otherwise not entirely clear
(we will give another more geometrical expression later). One may
further write it in terms of $\extd$ using results from the last
section. For a simpler expression, the second line in
Proposition~\ref{delbardel} gives
\[ qb_-\extd b_0\wedge \extd b_+-q^{-1}b_0\extd b_-\wedge
\extd b_+=q^2b_0^2\vol\] which gives the volume form if $b_0$ is
invertible. An alternative is to use (\ref{db}) and compute \[
\extd b_+\wedge \extd b_--q^4\extd b_-\wedge \extd
b_+=q^3[2]_q(1+[2]_q b_0)\vol\] which gives the volume form if one
assumes $1+[2]_qb_0$ invertible. These are classically the two
trivialisations of the sphere given by deleting the north or south
poles.

Next we look at the metric, motivated from Corollary~\ref{dSrel}.
We let $\bar\tens$ denote the tensor product over $\C_q[S^2]$.

\begin{propos}\label{g} There is a natural metric
\[ g=q^2\extd b_-\bar\tens \extd b_++\extd b_+\bar\tens \extd
b_--[2]_q\extd b_0\bar\tens \extd b_0\] such that $g$ is invariant
under the left coaction of $\C_q[SL_2]$ and $q$-symmetric in the
sense $\wedge(g)=0$. Moreover, $g\in (\Omega^{1,0}\bar\tens\Omega^{0,1})\oplus
(\Omega^{0,1}\bar\tens\Omega^{1,0})$. 
\end{propos}
\proof Left-invariance of the metric follows from
\[\begin{pmatrix}a^2& [2]_qab & b^2\\ ca & 1+[2]_qbc & db\\
c^2& [2]_qcd & d^2\end{pmatrix}^t\begin{pmatrix}0&0&q^2\\
0&-[2]_q&0\\ 1&0&0\end{pmatrix}  \begin{pmatrix}a^2& [2]_qab
& b^2\\ ca & 1+[2]_qbc & db\\
c^2& [2]_qcd & d^2\end{pmatrix}=\begin{pmatrix}0&0&q^2\\
0&-[2]_q&0\\ 1&0&0\end{pmatrix}\] where $t$ denotes transpose and
where we use the relations of $\C_q[SL_2]$. The transformation
matrix here is the one in the coaction (\ref{coactS}) in the basis
$\extd b_-,\extd b_0,\extd b_+$. Actually, this coaction
corresponds to the vector corepresentation of the even subalgebra
$\C_q[SO_3]$ of $\C_q[SL_2]$ and for generic values of $q$ there
is a unique invariant such matrix for the metric coefficients up
to a scale. Hence the metric is uniquely determined if we suppose
it has numerical coefficients with our basis of exact
differentials are viewed as spanning a 3-dimensional vector space
over $\C$ (invariance at this level then implies invariance when
viewed over $\bar\tens$). Such numerical coefficients in turn are
a natural assumption in view of (\ref{dSrel}). Differentiating
that, we see that $\wedge(g)=0$. 

Also, writing $g=g_{++}\oplus g_{+-}\oplus g_{-+}\oplus g_{--}$
for the decomposition according to Corollary~\ref{dhol}, we use the $q$-commutation
relations there and expressions for $\del b_i$ in (\ref{deld}) to compute
\begin{eqnarray*}
g_{++}\nquad &&=\del b_+\bar\tens \del b_--[2]_q\del b_0\bar\tens \del b_0+q^2\del b_-\bar\tens \del b_+\\
&&=q^3b_-\del b_+\bar\tens\extd b_0+(q^3-q^{-1})b_+\del b_-\bar\tens \extd b_0-q^3 b_0\del b_+\bar\tens \extd b_- +[2]_qb_+\del b_0\bar\tens \extd b_-\\
&&\quad -[2]_q(1+q^3b_0)\del b_0\bar\tens \extd b_0
+q^2(1+q b_0)\del b_-\bar\tens\extd b_+-q^{-1}b_+\del b_-\bar\tens \extd b_0\\
&&=-q^{-1}(1+q^2[2]_q b_0)\del b_0\bar\tens \extd b_0+(q^2[2]_q-q^3)b_0\del b_+\bar\tens\extd b_-+q^2(1+qb_0)\del b_-\bar\tens\extd b_+\\ &&=0
\end{eqnarray*}
on using $(1+q^2[2]_q b_0)\del b_0=(\del b_0)(1+[2]_qb_0)$ and then using Corollary~\ref{dSrel} to replace $\extd b_0$ by $\extd b_\pm$. Then, as well as for the third 
equality, we used
the relations in $\Omega^{1,0}$ in Corollary~\ref{dhol} to collect terms. There is
a similar proof for $g_{--}=0$. \eproof

Next, we look at the Hodge $*$ operator \[
*:\Omega^1(\C_q[S^2])\to \Omega^1(\C_q[S^2])\] which we require to
obey $*^2=\id$ and to be at least a left-module map and to be
frame-invariant (the metric can also be analysed in such terms but
frame invariance alone does not fix a particular one, we used
rotational left-covariance). In the frame bundle approach for $*$
we require $*: V\to V$ where $V$ is the 2-dimensional local
tangent space. In order to be frame invariant (which means
covariant under (\ref{coactV})) and square to the identity, this
must be given by $*(e^\pm)= \pm e^\pm$ up to an overall sign.

\begin{propos}\label{hodge} The natural Hodge $*$ operator is a
left-covariant bimodule map obeying
\[  *(\del f)= \del f,\quad *(\bdel f)=-\bdel f,\quad \forall
f\in \C_q[S^2].\]
and define a left-covariant lifting
$i:\Omega^2(\C_q[S^2])\to \Omega^1(\C_q[S^2])\bar\tens
\Omega^1(\C_q[S^2])$,
\[ i(\vol)={q^{-1}\over
[2]_q}(*\bar\tens\id)(g)= -{q^{-1}\over
[2]_q}(\id\bar\tens
* )(g)\]
\end{propos}
\proof Let us first verify directly that $*$ is well-defined as a
left module map, in which case it is given as stated since $\bdel
b_\pm,\bdel b_0$ are generated from $e^-$, etc. Indeed
$\Omega^1(\C_q[S^2])$ is a rank 2 bundle in which we can take
$\extd b_\pm, \extd b_0$ as generators with the relation in
Corollary~\ref{dSrel}. Writing $\extd=\del+\bdel$ we have each
part holding separately, \eqn{dShol}{ q^2b_-\bdel b_++b_+\bdel
b_--(1+[2]_qb_0)\bdel b_0=0} and similarly for $\del$ (this is
also clear from relations in the proof of Corollary~\ref{dhol}),
so $*$ is compatible with this relation. That $*$ is a right
$\C_q[S^2]$ bimodule map can easily be proven using the Leibniz
rule for $\extd,\bdel,\del$. Moreover, $*$ is left-covariant under
the coaction of $\C_q[SL_2]$ since the coaction acts on each
$\Omega^{0,1},\Omega^{1,0}$ separately. Next, we recall the usual
formulae in which the Hodge $*$ operator is given in terms of the
`totally antisymmetric tensor' and the metric. The role of that
tensor is played by the lifting of the volume form to an element
of $\Omega^1\bar\tens\Omega^1$ which is something that in
classical geometry one takes for granted (the wedge product is
given classically by skew-symmetrization). In noncommutative
geometry, as explained in \cite{Ma:rieq}, this lifting map is an
additional datum required to split the wedge map
$\Omega^1\bar\tens\Omega^1\to \Omega^2$ and we use the Hodge *
operator and the metric to define it. We check that
\begin{eqnarray*}&&\nquad \wedge(*\bar\tens\id)(g)
\\
&&=q^2(-a^2 e^-+b^2 e^+)\wedge (c^2 e^-+d^2 e^+)+(-c^2 e^-+d^2
e^+)\wedge (a^2 e^-+b^2 e^+)\\
&&\quad -[2]_q(-ca e^-+db e^+)\wedge (ca e^-+ db e^+)\\
&&=(q^2 a^2d^2 +q^4 b^2c^2+c^2b^2+ q^2 d^2 a^2 -[2]_q ca
db-q^2[2]_q db ca)\vol\\
&&=(q^2+1)\vol
\end{eqnarray*}
where we work `upstairs' in the frame bundle and where the even
terms give $q^2$ and the odd terms give $1$ using the relations of
$\C_q[SL_2]$. Hence $\wedge \circ i(\vol)=\vol$ as required. One
may also obtain this result using the computations in the proof of
Proposition~\ref{delbardel}. We define $i$ as by definition a
left module map (it is not a bimodule map). The other stated expression for $i$ 
is the same in view of the form of $g$ in Proposition~\ref{g}. 
Note also that since $\wedge(g)=0$ we have in fact a general family of
lifts of this type, \eqn{generali}{ i(\vol)={q^{-1}\over [2]_q}
(*\bar\tens\id)g+\mu g} for any $\mu$. Or equivalently, $i(\vol)=\alpha g_{+-}+\beta g_{-+}$ provided
$\alpha-\beta=2q^{-1}/[2]_q$. These lifts are all left covariant
under the coaction of $\C_q[SL_2]$ since $g_{+-},g_{-+}$ separately are. \eproof

We do not explicitly discuss complex structures in this paper; we
work over $\C$ but with care this could be any field as in
algebraic geometry. Nevertheless, our $b_\pm$ coordinates have
their interpretation as complex linear combinations (of the ambient
$\R^3$ coordinates) in real
geometry; in real coordinates the Hodge $*$ is equivalent to an
almost complex structure $J$ since this is defined in two
dimensions exactly by the same relation between the volume form
(viewed as a symplectic structure) and the metric as for $i$ 
in Proposition~\ref{hodge}. In our case since we are
deforming the standard metric on the sphere, this gives implicitly
a $q$-deformation of its actual complex structure via the Hodge *.
This justifies our notations $\del,\bdel$.

We now have all the basic structures at least for the first
`layer' of geometry, namely cohomology and electromagnetism.  For
the Maxwell theory we define of course
\[ *1=\vol, \quad *\vol=1\]
and the Laplacian on degree zero by $\square=-{1\over 2}*\extd
*\extd$. Then \eqn{lapsq}{ (\square f)\vol =-{1\over 2}\extd*\extd
f={1\over 2}\extd(\bdel f-\del f)=\del\bdel f.}

\begin{propos}\label{laplace} The functions $b_-,1+[2]_qb_0,b_+$ are
eigenfunctions of $\square$ with eigenvalue $q^2[2]_q$.
\end{propos}
\proof We compute
\[ \extd\del\begin{cases}b_+\\ b_-\\
b_0\end{cases}=\begin{cases} q\extd b_0\wedge \extd
b_+-q^{-1}\extd b_+\wedge\extd b_0\\
q\extd b_-\wedge \extd b_0-q^{-1}\extd b_0\wedge \extd b_-\\
q^2\extd b_-\wedge\extd b_+-q^{-1}\extd b_0\wedge\extd
b_0\end{cases}\] since $\extd=\bdel$ on $\Omega^{1,0}$.  We
compute these using (\ref{db}) in the same manner as in the proof
of Proposition~\ref{hodge}. Here we need their values explicitly:
\begin{eqnarray*} \extd b_-\wedge \extd b_+&=&-(1+q^{-2}[2]_q
b_0-q^{-1}(q-q^{-1})[3]_q b_0^2)\vol\\
\extd b_0\wedge\extd b_0&=&(q-q^{-1})q^2([2]_q+[3]_qb_0)b_0\vol\\
\extd b_0\wedge\extd b_+&=&((q^5-q^{-1})b_0-1)b_+\vol\\
\extd b_+\wedge\extd b_0&=&( (q^7-q)b_0+q^4)b_+\vol\\
\extd b_-\wedge\extd b_0&=&b_-((q^5-q^{-1})b_0-1)\vol\\
\extd b_0\wedge \extd b_-&=&(q^7-q)b_- b_0\vol
\end{eqnarray*}
which we use to find $\extd\del=\bdel\del=-\del\bdel$. \eproof

Note that these functions form the vector corepresentation of
$\C_q[SO_3]$ under the left coaction (\ref{coactS}) (this appears
in the above basis as the transformation matrix in the proof of
Proposition~\ref{g}). In the same way, the matrix elements of each
integer spin corepresentation $V_n$ of $\C_q[SL_2]$ define a
square-dimension subspace $V_n\tens V_n^*\subset\C_q[SL_2]$.
Fixing the unique zero weight vector $v$ under the right coaction,
the subspace $V\tens v$ (in other words, the matrix entries in the
middle row of the transformation matrix) span an eigenspace of the
Laplacian. For generic $q$ the Peter-Weyl decomposition of
$\C_q[SL_2]$ implies, as classically, that this is a complete
diagonalisation of $\square$ on $\C_q[S^2]$ with one eigenspace
for each integer spin corepresentation. The matrix elements of the
$1/2$-integer spin corepresentations contain an odd number of the
$\C_q[SL_2]$ generators and hence can never have the zero degree
needed to lie in $\C_q[S^2]$.

In particular, the zero eigenspace of $\square$ is spanned by the
constant function $1$. This implies that if $\bdel f=0$ then $f$
is a multiple of $1$. Similarly for $\del$, so for the
noncommutative de Rham and Dolbeault cohomology of $\C_q[S^2]$ for
generic $q$ we have
\[ H^0=H_\del^0=H^0_{\bdel}=\C.1.\]
It is clear that for generic $q$ we also have the usual values for
the rest of the cohomology as well as Poincar\'e duality, since
all constructions for this calculus are a smooth deformation of
the classical ones.  We omit explicit proofs of these facts since
we need them for discussion only.

Similarly, we have a ``massive'' Maxwell equation defined on
1-forms by
\[ \square_1A\equiv -{1\over 4}*\extd*\extd A= m^2 A,\quad *\extd * A=0\]
The second equation is Coulomb gauge in physics and is automatic
when $m\ne 0$ (in this case $m^2A$ should be interpreted as the
source). We recall that in Maxwell theory the field is considered
modulo exact forms but this freedom can be partially fixed by a
gauge choice. We write the curvature as $F=\extd A=f\vol$ where
$f=*F$ is in $\C_q[S^2]$. Then $\square_1A=0$ translates to
\[ \del f=\bdel f=0\]
which implies that $f\propto 1$. In that case, $\extd A\propto
\vol$ which implies $f=0$ since $\vol$ is not exact by Poincar\'e
duality, so the only `photons' are pure gauge. This is to be
expected for a sphere. On the other hand, if $A$ is a "massive"
mode, then
\[ A=-{1\over m^2}*\extd f={1\over 4m}(\bdel f-\del f),\quad
\square f=2m^2 f.\] Conversely, given an eigenfunction of
$\square$ as in the second equation, we use the first to define
$A$ and obtain a massive mode. For example, the eigenfunctions of
$\square$ in Proposition~\ref{laplace} give solutions
\[ A={q^2\over 2[2]_q}(\bdel b_i-\del b_i)\]
where $i=\pm,0$ and these are given as 1-forms via (\ref{deld}).

\section{Levi-Civita connection, curvature and Dirac
operator on the $q$-sphere}

Next, we compute use the frame bundle approach to develop the
Riemannian geometry of the $q$-sphere. Here the 4$q$-monopole
bundle \cite{BrzMa:gau} is viewed as the frame bundle and the
$q$-monopole connection (\ref{omegamon}) on it as a spin
connection. We find that it correctly induces the Levi-Civita
connection on the cotangent bundle.

\begin{theorem}\label{cov} The $q$-monopole connection
(\ref{omegamon}) viewed as spin connection in the frame bundle of
$\C_q[S^2]$ induces  the covariant derivative
\[ \nabla \begin{cases}\extd b_\pm\\
\extd b_0\end{cases}= \left\lbrace\begin{array}{c} [2]_qb_\pm\\
1+[2]_q b_0\end{array}\right\rbrace g\]
which is torsion-free and skew-metric compatible in the sense of
zero cotorsion (a generalised Levi-Civita connection). 
\end{theorem}
\proof We recall that the bundle $\CE_{-2}=(\C_q[SL_2]\tens
[b_-])^{\C[t,t^{-1}]}$ can be identified with $\C_q[SL_2]_2\tens
[b_-]$, where we now write the representative $[b_-]\in V$
explicitly. This space in turn was identified with
$\C_q[SL_2]_2.e^-=\Omega^{0,1}$ as explained in the proof of
Corollary~\ref{dhol}. The action of the covariant derivative on
$\CE_{-2}$ is by \[ D(f\tens [b_-])=(\id-\Pi_\omega)(\extd f)\tens
[b_-]=(\extd f-f\omega(t^2))\tens [b_-]=(\extd f-(1+q^2)fe^0)\tens
[b_-]\] for all $f\in \C_q[SL_2]_2$. This is the usual covariant
derivative on a $q$-monopole section, here of charge 2. Working
`upstairs' on $\C_q[SL_2]$ and using the Leibniz rule and the 3-d
calculus, we have
\[ \extd (a^2)=(1+q^2)(a^2 e^0+qab e^+),\quad \extd(ca)=(1+q^2)ca
e^0+q^2(1+[2]_qbc)e^+\]
\[\extd(c^2)=(1+q^2)(c^2 e^0+qcd e^+).\]
We see that the horizontal projection simply kills the $e^0$ term
in each expression. Next, from the identity
$d^2a^2+q^2b^2c^2-q[2]_qdbac=1$ we write
\eqn{eplus}{e^+=1.e^+=q^{-2}\del b_+.a^2+\del
b_-.c^2-q^{-2}[2]_q\del b_0.ca} where we move the $a,c$ generators
to the right using the relations of the 3-d calculus. These degree
2 products combine with $[b_-]$ to give a section of $\CE_{-2}$.
We are working `upstairs' but we can now identify the product as
the $\C_q[S^2]$-module structure. Thus \[ D(a^2\tens
[b_-])=[2]_qb_-\left(\del b_+.(a^2\tens [b_-])+q^2\del
b_-.(c^2\tens [b_-])-[2]_q\del b_0.(ca\tens [b_-])\right).\]
Similarly for $ca\tens [b_-]$ and $c^2\tens [b_-]$. Finally, we
replace $[b_-]$ by $e^-$ (the framing isomorphism
Theorem~\ref{Sman}) and identify the resulting elements of
$\Omega^{0,1}$ on the right. This gives $\nabla \bdel b_-$.
Similarly for all the other cases. For the $\nabla\del$ we
use
\[ e^-=q^2\bdel b_-.d^2+\bdel b_+.b^2-[2]_q\bdel b_0.db.\]
As a result, we find
\eqn{covg}{  \nabla(\del b_\pm)=[2]_q b_\pm
g_{-+},\quad \nabla(\bdel b_\pm)=[2]_q b_\pm g_{+-}}
where 
\[ g_{+-}=q^2\del b_-\bar\tens\bar\del b_++\del b_+\bar\tens \bar\del b_--[2]_q\del b_0
\bar\tens\bar\del b_0\]
 etc. in the decomposition of $g$ as in the proof of Proposition~\ref{g}. Combining these gives $\nabla \extd$ as stated. 

Next, having found $\nabla$, we look at the torsion equation. As
explained in \cite{Ma:rie} the noncommutative meaning of this is
\[ {\rm Tor}=\nabla_\wedge - \extd: \Omega^1\to \Omega^2\] which is
the first degree measure of the failure of $\nabla\wedge$ to form
a complex (the second degree measure is the curvature). We have
\[{\rm Tor}(\extd b_\pm)=\nabla\wedge(\extd
b_\pm)=[2]_q b_\pm \wedge(g)=0\]
by the $q$-symmetry in Proposition~\ref{g}. Since the torsion is a left module
map, it follows that the torsion vanishes entirely.

Finally, we look at the `skew-metric compatibility' in the sense
of zero cotorsion. This has been proposed as the correct notion of
compatibility in \cite{Ma:rie} and can be written in terms of
\[ {\rm CoTor}=(\nabla\wedge\id-\id\wedge\nabla)g\in
\Omega^2\bar\tens\Omega^1\]
(there is an additional term if the torsion is not zero). Since
the metric consists of exact differentials and since the torsion
vanishes, the first $\nabla\wedge\id$ vanishes. Looking at the
second term, we compute
\begin{eqnarray*}{1\over [2]_q}(\id\wedge\nabla)g\nquad
&&={1\over[2]_q}(q^2\extd b_-\wedge \nabla\bar\tens \extd
b_++\extd b_+\wedge\nabla\bar\tens \extd b_-)-\extd
b_0\wedge\nabla\bar\tens
\extd b_0\\
&&=q^2\extd b_-\wedge b_+ g+\extd b_+\wedge b_-g-\extd b_0\wedge (1+[2]_qb_0)g\\
&&=0
\end{eqnarray*}
by a right-module version $q^2 (\extd b_-) b_++(\extd b_+) b_--\extd b_0(1+[2]_qb_0)=0$ 
of the relation in Corollary~\ref{dSrel}. Hence the cotorsion vanishes as well.
\eproof

Note that the cotorsion or `skew-metric compatibility' condition
appropriate in noncommutative geometry\cite{Ma:rie} is weaker than
the usual notion. In our case we have the more usual $\nabla g=
O(q-1)$ (if $\nabla$ is taken to act on the tensor product as a
derivation while keeping its left output to the far left), so that
we recover the usual full metric compatibility only when $q=1$. It
is also worth noting that viewed as sections of an associated
bundle (see the Appendix), it is the $\bdel b_\pm,\bdel b_0$ which
are actually holomorphic in the sense $\nabla \del b_i\in
\Omega^{1,0}\bar\tens \Omega^1$, rather than the image of $\del$.
Let us also use the connection to relate to the projective module
point of view on quantum bundles.
\begin{corol}
The projector
\[ E=\begin{pmatrix}[2]_q b_- \\ 1+[2]_q b_0 \\ [2]_q b_+ \end{pmatrix}\begin{pmatrix}-b_+,\ 1+[2]_q b_0,\  -q^2 b_-\end{pmatrix}\]
yields $\Omega^1(\C_q[S^2])=\C_q[S^2]^{\oplus 3}.(1-E)$ and $\nabla=-E\extd E$ acting on $\extd b_-,\extd b_0,\extd b_+$.
\end{corol}
\proof Note that proceeding from Theorem~\ref{Sman} would give a projector
from 6 copies of $\C_q[S^2]$ whereas we provide a projector more in keeping
with the classical geometrical picture from 3 copies. Moreover, we are not using the universal calculus as in \cite{HajMa:pro}. Nevertheless, the form of $\nabla$ similarly suggests the projection shown, which we then verify directly. Thus, we have the dot products
\[ \begin{pmatrix}-b_+,\ 1+[2]_q b_0,\ -q^2 b_-\end{pmatrix}\begin{pmatrix}[2]_q b_- \\ 1+[2]_q b_0 \\ [2]_q b_+ \end{pmatrix}=1,\quad \begin{pmatrix}-b_+,\ 1+[2]_q b_0,\ -q^2 b_-\end{pmatrix} \begin{pmatrix}\extd b_- \\ \extd b_0 \\ \extd b_+ \end{pmatrix}   =0\]using respectively, the relations (\ref{Srel}) of the $q$-sphere and the relation in Corollary~\ref{dSrel}. The second dot product with $\extd$ of the row vector similarly gives $-g$. These observations imply that $E^2=E$ and (using the Leibniz rule to compute $\extd E$) that $\nabla =-\extd E=-E\extd E $ when acting on the column vector $(\extd b_i)$. Here
$E.(\extd b_i)=0$ (acting on the column vector). The map $\C_q[S^2]^{\oplus 3}\to \Omega^1(\C_q[S^2])$ is given by $(f,g,h)\mapsto f\extd b_-+g\extd b_0+h\extd b_+=(f,g,h)(1-E)(\extd b_i)$ and has kernel generated as a left module by the row vector in $E$. This corresponds to the relation in Corollary~\ref{dSrel}. Let us note that the same
$E$ also gives the $\del$ and $\bar\del$ parts in a similar way. Thus
$\nabla=-\del E=-E\del E$ when acting on the column vector $(\bar\del b_i)$ and $\nabla=-\bar\del E=-E\bar\del E$ on $(\del b_i)$. One may check that $\bar\del E.(\bar\del b_i)=\del E.(\del b_i)=0$ so that $\nabla=-\extd E=-E\extd E$ when acting on either $(\del b_i)$ or $(\bar\del b_i)$ separately. \eproof

\begin{propos}\label{riemann} The Riemann and Ricci tensors of the
above generalised Levi-Civita connection are \[   {\rm
Riemann}|_{\Omega^{0,1}}=[2]_q\vol\bar\tens \id,\quad {\rm
Riemann}|_{\Omega^{1,0}}=-q^4[2]_q\vol\bar\tens\id.\] The lift
\[ i(\vol)={q^{-1}\over [2]_q}\left(-(\id\bar\tens
*)g+{1-q^{-4}\over 1+q^{-4}}g\right)\] and trace in the middle
position gives
\[ {\rm Ricci}={2 q^{-1}\over 1+q^{-4}}g\]
making the $q$-sphere an `Einstein space'.
\end{propos}
\proof The Riemann tensor is defined
abstractly\cite{Ma:rie,Ma:rieq} by
\[ {\rm Riemann}=(\id\wedge\nabla-\extd\bar\tens\id)\nabla:\Omega^1\to
\Omega^2\bar\tens\Omega^1\] as the form-version of the usual
definition, as explained in \cite{Ma:rie}. One may compute it from
the formulae for $\nabla$ above. Since $\nabla(\bdel b_\pm)\in
\Omega^{1,0}\bar\tens\bdel\{b_i\}$, when we apply
$\id\wedge\nabla$ we will get zero since $\Omega^{2,0}=0$. So only
the $-\extd\tens \id)\nabla$ contributes. We have
\begin{eqnarray*}-{1\over [2]_q}(\extd\bar\tens\id)\nabla(\bdel
b_+)\nquad&&=-b_+(q^2\bdel\del b_-\bar\tens \bdel b_++\bdel\del
b_+\bar\tens\bdel b_--[2]_q\bdel\del b_0\bar\tens\bdel b_0)\\
&&\quad -\bdel b_+\wedge(q^2\del b_-\bar\tens \bdel b_++\del
b_+\bar\tens \bdel b_--[2]_q\del b_0\bar\tens\bdel b_0).
\end{eqnarray*}
We use Proposition~\ref{laplace} for the Laplacian $\bdel\del$ and
that $\vol$ is central in the first line to collect $q^2[2]_q
b_+\vol$ to the left times (\ref{dShol}), so that the first line
vanishes. For the second line we use our computations for such
wedge products in  the proof of Proposition~\ref{delbardel} as
multiples of $\vol$, to obtain
\[-{1\over [2]_q}(\extd\bar\tens\id)\nabla(\bdel
b_+)=\vol (q^2b_0^2\bar\tens\bdel b_++q^{-1}b_+^2\bar\tens\bdel
b_--[2]_qq^{-1}b_+b_-\bar\tens\bdel b_0=\vol\bar\tens \bdel b_+\]
on using the relations of the $q$-sphere and the relations between
the $\bdel b_\pm,\bdel b_0$ in Corollary~\ref{dhol}. Similarly for
the Riemann tensor on $\bdel b_-$. The computation for ${\rm
Riemann}(\del_\pm)$ is similar but yields an extra factor $-q^4$
(the symmetry was broken in our choice of $\vol$). We note that
Riemann is a left module map so it is enough to find it on such
exact differentials. It is also possible to compute the curvature
`upstairs' in the principal bundle, using (\ref{Fmon}). By the
same conventions and explanations as in the proof of
Theorem~\ref{cov}, we have, for example \[{\rm Riemann}(\del
b_+)=b^2F_\omega(t^{-2}) e^+=q^3[-2;q^2]b^2\vol\bar\tens
e^+=-q^4[2]_q\vol \bar\tens \del b_+.\] The curvature can be
computed either way, as explained for the classical case in
\cite{Ma:rie}. Using the Hodge $*$ operator we can write the
Riemann tensor as
\[ {\rm Riemann}=[2]_q\vol\bar\tens({1-q^4\over 2}-{1+q^4\over 2}*).\]

For the Ricci tensor we need to lift the Riemann tensor to a map
$\Omega^1\to \Omega^1\bar\tens\Omega^1\bar\tens\Omega^1$,
\begin{eqnarray*} i({\rm Riemann})&=&[2]_q i(\vol)\bar\tens ({1-q^4\over
2}-{1+q^4\over 2}*)\\
&=&{2 q^{-1}\over 1+q^{-4}}\left(\id\bar\tens ( {1-q^{-4}\over
2}-{1+q^{-4}\over 2}*)\right)g \bar\tens ({1-q^4\over
2}-{1+q^4\over 2}*)\end{eqnarray*} where we use the lifting from
the same family as in Proposition~\ref{hodge} but of the form
stated. We can then take a trace by `feeding' the right hand
factor of $i(\vol)$ into the input of Riemann.  This gives
\[ {\rm Ricci}= {2 q^{-1}\over 1+q^{-4}}\left(\id\bar\tens
({1-q^4\over 2}-{1+q^4\over 2}*)( {1-q^{-4}\over
2}-{1+q^{-4}\over 2}*)\right)g
\] which gives the result stated using $*^2=\id$. \eproof

The lift $i:\Omega^2\to \Omega^1\bar\tens\Omega^1$ is an
additional datum in the approach of \cite{Ma:rieq} needed to
define the Ricci tensor as well as interior products etc. Our
point of view in the above Proposition~\ref{riemann} is that for
the standard metric to be Einstein, the natural lift $i$ in
Proposition~\ref{hodge} gets deformed by an additional $g$
component which vanishes as $q\to 1$. Equivalently, if we keep the 
choice of $i$ coming from the geometry in Proposition~\ref{hodge} 
then we find \eqn{altricci}{
i(\vol)=-{q^{-1}\over [2]_q}(\id\bar\tens *)g \quad
\Rightarrow\quad {\rm Ricci}={q^{-1}(1+q^4)\over
2}g+{[2]_q(1-q^4)\over 2}i(\vol)} which is a novel prediction of
an `antisymmetric' volume form correction to the Ricci tensor that
vanishes as $q\to 1$.

Let us also note that the definition of Ricci used above is via
the trace as in \cite{Ma:rieq}, but between the second factor of
the lift of Riemann (rather than the first factor as there) and
its input. Also, we do not need a
braided trace as was needed to keep covariance in the bicovariant
calculus model in \cite{Ma:ric}, and do not have offsets
$\theta\tens\theta$ in Ricci as appeared there.   The completely
general definition of Ricci at the level of arbitrary framed
algebras is not fully understood, but we see once again that in
examples, as here, it is clear which trace to take.

The $\gamma$ matrices needed for the Dirac operator are likewise
not yet formulated in the most general form for any framed
algebra, but in examples there seems to be a clear choice. We
propose the following. For the spin bundle we take \[
\CS\equiv\CS_-\oplus\CS_+=\CE_{-1}\oplus
\CE_{+1}=\C_q[SL_2]_{1}\oplus\C_q[SL_2]_{-1}\] as given by the the
monopole bundles of charges -1 and 1. We identify the sections of
the bundles with the degree $\pm1$ subspaces of $\C_q[SL_2]$ as we
did already for the charges $\pm2$. This corresponds to the double
cover of the bundle for the cotangent space in Theorem~\ref{Sman}.
The next ingredient
 is a map
\[ \gamma:\Omega^1(\C_q[S^2])\to {\rm End}(\CS)\]
which we construct as follows. We use our description in
Corollary~\ref{dhol} as we have throughout the paper to define

\eqn{gammadef1}{\gamma:\Omega^{1,0}\bar\tens\CS_-\to \CS_+,\quad
\gamma(fe^+\bar\tens \sigma)=f\sigma,\quad \forall f\in
\C_q[SL_2]_{-2},\quad \sigma\in \CS_-}\eqn{gammadef2}{
\gamma:\Omega^{0,1}\bar\tens\CS_+\to \CS_-,\quad
\gamma(fe^-\bar\tens \tau)=f\tau,\quad \forall f\in
\C_q[SL_2]_2,\quad \tau\in \CS_+.} Here $\sigma,\tau$ denote
appropriate sections and $\gamma$ under our identifications is
nothing other than the product of $\C_q[SL_2]$ restricted to the
appropriate degrees. We also let \eqn{gammadef3}{
\gamma|_{\Omega^{0,1}\bar\tens \CS_-}=0,\quad
\gamma|_{\Omega^{1,0}\bar\tens \CS_+}=0.} The classical motivation
for $\gamma$ is as follows. Since $\del$ is like a holomorphic
differential one may think of it is a complex linear combination
of the usual differentials. Likewise, if $\sigma^{1,2}$ are the
usual Pauli matrices, then \[ {1\over
2}(\sigma^1+\imath\sigma^2)=\begin{pmatrix}0&1\\
0&0\end{pmatrix},\quad {1\over
2}(\sigma^1-\imath\sigma^2)=\begin{pmatrix}0&0\\
1&0\end{pmatrix}\] which is the structure we have used for
(\ref{gammadef1})--(\ref{gammadef3}).

\begin{lemma}\label{gamma}
The  operator $\gamma:\Omega^1(\C_q[S^2])\bar\tens \CS\to \CS$
defined above is covariant under the left coaction of $\C_q[SL_2]$
and obeys
\[ \{\gamma_{\extd b_\pm},\gamma_{*\extd b_\pm}\}=0,\quad
\gamma_{\extd b_\pm}\circ\gamma_{\extd b_\pm}=\begin{pmatrix}q^{-1}&0\\
0 & q^{3}\end{pmatrix}b_\pm^2=-\gamma_{*\extd
b_\pm}\circ\gamma_{*\extd b_\pm} \] where $\gamma_{\extd b
_\pm}=\gamma(\extd b_\pm\bar\tens(\ ))$, etc. Moreover,
\[ \gamma\circ\gamma(g\bar\tens(\ ))=\begin{pmatrix} q^2 & 0
\\ 0 & 1\end{pmatrix} \id\]
\end{lemma}
\proof The left coaction on our various spaces is simply the
coproduct $\Delta$ of $\C_q[SL_2]$ restricted to the appropriate
degree. Since this is an algebra homomorphism, the $\gamma$ as
given by the product is covariant. Likewise, $\del b_i,\bdel b_i$
correspond as in Corollary~\ref{dhol} to $b^2,db,d^2,a^2,ca,c^2$
and the relations among the corresponding $\gamma_{\del
b_i},\gamma_{\bdel b_i}$ are just the relations among these
generators of $\C_q[SO_3]$ defined as the even part of
$\C_q[SL_2]$, since $\gamma$ acts by left multiplication. Thus
\begin{eqnarray*} \gamma_{\del b_-}\circ\gamma_{\bdel b_-}
=b^2a^2=q^3 b_-^2,&\quad &\gamma_{\bdel b_-}\circ\gamma_{\del b_-}
=q^{-1}b_-^2\\
\gamma_{\del b_+}\circ\gamma_{\bdel b_+}=d^2c^2=q^3 b_+^2,
&\quad &\gamma_{\bdel b_+}\circ\gamma_{\del b_+}=q^{-1}b_+^2\\
\gamma_{\del b_0}\circ\gamma_{\bdel
b_0}=dbca=q^2(1+qb_0)b_0,&\quad
&\gamma_{\bdel b_0}\circ\gamma_{\del b_0}=(1+q^{-1}b_0)b_0 \\
\gamma_{\del b_+}\circ\gamma_{\bdel b_-}=d^2a^2=
(1+q^3b_0)(1+qb_0),&\quad &\gamma_{\bdel b_-}\circ\gamma_{\del
b_+}=(1+q^{-3}b_0)(1+q^{-1}b_0)\\
\gamma_{\del b_-}\circ\gamma_{\bdel b_+}=b^2c^2=b_0^2,&\quad
&\gamma_{\bdel b_+}\circ\gamma_{\del b_-}=b_0^2.
\end{eqnarray*}
The left column act on $\CS_+$ and the right column on $\CS_-$ by
multiplication. Remembering that $\gamma$ acts by zero when the
degrees do not match, we find $\gamma_{\extd b_\pm}=\gamma_{\del
b_\pm}+\gamma_{\bdel b_\pm}$ with square as stated. Since the
Hodge
* changes the sign of $\bdel b_i$ we find similarly that
$\gamma_{\extd b_\pm}$ and $\gamma_{*\extd b_\pm}$ anticommute.
Finally, the other expressions allow us to compute
\[ \gamma\circ\gamma( g_{+-}\bar\tens(\ ))=q^2b_0^2+(1+q^3b_0)(1+q
b_0)-[2]_qq^2(1+q b_0)b_0=1,\quad
\gamma\circ\gamma(g_{-+}\bar\tens(\ ))=q^2\] where for example
$g_{+-}$ is the $\Omega^{1,0}\bar\tens\Omega^{0,1}$ component of
the metric. These combine to the result stated. \eproof

The classical meaning of the $\gamma$ relations stated is that in
a local coordinate chart the 2-dimensional cotangent space is
spanned by $\extd b_+,*\extd b_+$, say, or (another chart) with
$b_-$. We see that our $\gamma$ operators in these directions
mutually anticommute and each square to a multiple of a
$q$-deformation of the identity. The relation involving the metric
is a weak form of the Clifford relations involving the metric, as
proposed in \cite{Ma:rieq}.

\begin{propos}\label{Dsl} $\phantom{\vardelta}$ Let $D$ be the
covariant derivative on $\CS$ given by the $q$-monopole as spin
connection. We define the Dirac operator on $\C_q[S^2]$ by
\[ \Dsl=\gamma\circ D=\begin{pmatrix}0& \gamma\circ\bar\vardelta \\
\gamma\circ\vardelta & 0\end{pmatrix}:\CS\to \CS\] where
$D=\vardelta+\bar\vardelta$ according to the parts in
$\Omega^{1,0}$ and $\Omega^{0,1}$. Then $\Dsl$ is covariant under
the left coaction of $\C_q[SL_2]$ and under local frame rotations
$\C[t,t^{-1}]$. Moreover, for $f=b_-,1+[2]_qb_0,b_+$, we have
\[ \Dsl^2(fa)=q^{-1}[2]_q(\square f)a+\begin{cases}0\\-q^{-1}a\\
-q^{-1}c,\end{cases}
\quad \Dsl^2(fb)=q^{-1}[2]_q(\square f)b+\begin{cases}0\\qb\\
qd\end{cases}\]
\[ \Dsl^2(fc)=q^{-1}[2]_q(\square f)c+\begin{cases}a\\ qc\\ 0,\end{cases}
\quad \Dsl^2(fd)=q^{-1}[2]_q(\square f)d+\begin{cases}-q^2b\\-q^3d\\
0.\end{cases}\]
\end{propos}
\proof If $\sigma\in \CS_-$ then $\gamma(D\sigma)=\gamma(\vardelta
\sigma)\in \CS_+$ etc., giving the stated form of $\Dsl$ on
$\CS_-\oplus\CS_+$. The space $\CS_-$ viewed as the degree 1
subspace of $\C_q[SL_2]$ is spanned over $\C_q[S^2]$ by $a,c$.
These are not linearly independent but obey the relations \[ b_+
a-(1+q b_0) c=0,\quad b_0 a-q^2 b_-c=0\] as used in the projector
\cite{HajMa:pro}. The covariant derivative on such sections is
known already from \cite{BrzMa:gau} and takes the form \[ Da=\extd
a-ae^0=qbe^+=q^{-1}\del b_0.a-q\del b_-.c ,\quad Dc=\extd
c-ce^0=qde^+=\del b_+.a-q\del b_0.c\] by similar computations as
in Theorem~\ref{cov}. We omit writing the basis of the degree -1
left comodule $V=\C$ in view of our identifications. Then
\[ \Dsl a=\gamma(bd e^+ \bar\tens a)-q\gamma(b^2 e^+\bar\tens
c)=bda-qb^2c=b.\] By such calculations, one has \eqn{Dslabcd}{\Dsl
a=b,\quad \Dsl c=d,\quad \Dsl b=qa,\quad \Dsl d=qc} where $b,d\in
\CS_+$ and $a,c\in \CS_-$.

For $\Dsl^2$ we note the Leibniz rule for any $\sigma\in \CS_-$
(say) and $f\in \C_q[S^2]$, \eqn{Dslleib}{ \Dsl(f\sigma)=f
\Dsl\sigma +\gamma(\del f\bar\tens\sigma)=f\Dsl\sigma +
f_i\gamma(\del b_i\bar\tens\sigma)} where $\del f=f_i \del b_i$
(sum over $i=-,0,+$) say. We similarly write $\bdel f= f_i
\bdel b_i$ and have a similar expression for the Liebniz property
on $\CS_+$. We choose the coefficients $f_i$ 
from a fixed expansion $\extd f=f_i\extd b_i$ (they are not unique). Then
\[ \Dsl^2(f\sigma)=f\Dsl^2\sigma+f_i\gamma(\bdel b_i\bar\tens
\Dsl\sigma)+\gamma(\bdel f_i\bar\tens\gamma(\del b_i\bar\tens
\sigma))+ f_i\Dsl\circ\gamma(\del b_i\bar\tens \sigma) .\] From
Theorem~\ref{cov} we have
\[ \nabla\del f=\bdel f_i\bar\tens \del b_i+
f_i\nabla\del b_i=\bdel f_i\bar\tens \del b_i+ f_i[2]_q b_i
g_{-+}+ f_0 g_{-+}\] which combined with Lemma~\ref{gamma} gives
the third term of $\Dsl^2$ as
\[\gamma\circ\gamma(\nabla\del f\bar\tens\sigma)-(q^2
f_i[2]_q b_i+q^2 f_0)\sigma.\] Meanwhile direct computation gives
\[ \Dsl\circ\gamma(\del
b_i\bar\tens\sigma)=\Dsl\left(\begin{cases}b^2\\ db\\
d^2\end{cases}.\sigma\right)=(q^2[2]_q b_i +q^2
\delta_{i,0})\sigma\] at least when $\sigma=a,c$. Hence
\[ \Dsl^2 (fa)=qfa+f_i\gamma(\bdel b_i\bar\tens
b)+\gamma\circ\gamma(\nabla\del f\bar\tens a)\]   and similarly
for $\Dsl^2(fc)$. Computing the middle terms, we find
\begin{eqnarray} \Dsl^2 (fa)&=&qfa+q^{-1}f_i b_i a-q^{-1}f_+c
+\gamma\circ\gamma(\nabla\del f\bar\tens a)\nonumber\\
\label{Dslsqac} \Dsl^2 (fc)&=& qfc+q^{-1}f_i b_i c+
f_-a+f_0c +\gamma\circ\gamma(\nabla\del f\bar\tens
c)\end{eqnarray} for all $f\in \C_q[S^2]$. There are similar
formulae for $\Dsl^2(fb),\Dsl^2(fd)$. The particular cases of $f$
stated then follow using Theorem~\ref{cov} or (\ref{covg}) and
$\gamma\circ\gamma(g\bar\tens(\ ))$ from Lemma~\ref{gamma}. We
recognise the action of the Laplacian from
Proposition~\ref{laplace}. One may also obtain these particular
cases by rather tedious direct computation, first finding \[
\Dsl(b_i a)=q[2]_q b_i b+\begin{cases}0\\ qb\\ d,\end{cases}\quad
\Dsl(b_i b)=[2]_q b_i a+\begin{cases}0\\ 0\\ -q^{-1}c\end{cases}
\]
\[ \Dsl(b_i
c)=q[2]_q b_i d+\begin{cases}-qb\\ 0\\ 0,\end{cases}\quad \Dsl(b_i
d)=[2]_q b_i c+\begin{cases}a\\ c\\ 0.\end{cases}
\]
where the cases are as $b_i=b_-,b_0,b_+$.
 \eproof

The formulae for $\Dsl^2$ are examples of  `Lichnerowicz type'
formulae, where relative to a chosen (holomorphic or
antiholomorphic) section, it is given by the scalar Laplacian from
Proposition~\ref{laplace} plus an additional term, which we think
of as some form of `scaler curvature'. Let us also comment that
our Dirac operator comes with a $\Z_2$ grading in the splitting
$\CS=\CS_-\oplus\CS_+$ and $\Dsl$ anticommutes with the grading.
Also, from (\ref{Dslleib}) and its cousin on $\CS_+$, we have
\eqn{Dslcom}{ [\Dsl,\hat f]=\gamma_{\extd f},\quad \forall f\in
\C_q[S^2]} where $\hat f$ denotes $f$ acting on $\CS$ by left
multiplication, so that $\Dsl$ does allow one to recover the
$\extd$ of the calculus at the algebraic level. Thus we have some
of the features of a spectral triple\cite{Con} although not
fitting precisely into that setting. On the other hand, our $\Dsl$
naturally deforms the geometrical Dirac operator on the classical
sphere and is motivated in that way rather than by such axioms.
Also, it is not hard to exhibit some eigenfunctions of $\Dsl$.
From the proof of Proposition~\ref{Dsl}, we have \eqn{Dsleval}{
\Dsl\begin{pmatrix}\pm q^\h a\\ b \end{pmatrix}=\pm
q^\h\begin{pmatrix} \pm q^{\h}a\\ b\end{pmatrix},\quad
\Dsl\begin{pmatrix} \pm q^{\h}c\\ d\end{pmatrix}=\pm
q^\h\begin{pmatrix}\pm q^{\h}c\\ d \end{pmatrix}} as solutions of
the massive Dirac equation.

Finally, let us note that unlike the cotangent bundle, the spinor bundle
is trivial. Both are trivial in $K$-theory having zero total monopole
charge (see \cite{HajMa:pro}).

\begin{propos}
Let $e=\begin{pmatrix}-q^{-1}b_0 & qb_- \\  -b_+ & 1+qb_0\end{pmatrix}$. Then
\[ \CS_-\isom \C_q[S^2]^{\oplus 2}(1-e),\quad \CS_+\isom \C_q[S^2]^{\oplus 2}e,\quad
\CS\isom \C_q[S^2]^{\oplus 2}.\]
The covariant derivative and Dirac operators in the trivialisation are
\[ D=\extd +(\extd e)e-e\extd e.\]
\begin{eqnarray*}\Dsl\nquad &&=\lambda q +\lambda q^{-1}(\lvec \del_i b_i)(1+(q^4-1)e)
\\
&& \quad +\lambda q^2\begin{pmatrix}\lvec\del_0& q^{-1}\lvec\del_+\\ -\lvec\del_-&0\end{pmatrix}e-\lambda \begin{pmatrix}0& q^{-1}\lvec \del_+\\ 
-\lvec\del_-&- \lvec\del_0\end{pmatrix}(1-e)\end{eqnarray*}
acting on row 2-vectors, where $\extd f=(f\lvec \del_i)\extd b_i$ (sum $i=-,0,+$) and $\lambda=q^{-\h}$.
\end{propos}
\proof  The projection for one half, the bundle $\CS_+$, say, was obtained
in \cite{HajMa:pro} with the
universal calculus and this is our starting point. Our observation
is that the projector for the other half of the spinor space
 is just given by the complementary projection. We then verify the 
desired properties directly. As spanning set for $\CS_-,\CS_+$ we take the
column vectors $\begin{pmatrix}a\\ c\end{pmatrix}$ and $\begin{pmatrix}b\\ d\end{pmatrix}$ respectively (these should not be confused with the $\CS_-\oplus\CS_+$ column vectors
above). We verify that
\[ e\begin{pmatrix}a\\ c\end{pmatrix}=0,\quad 
D\begin{pmatrix}a\\ c\end{pmatrix}=-\del e\begin{pmatrix}a\\ c\end{pmatrix}=-e\extd e \begin{pmatrix}a\\ c\end{pmatrix}\]
\[ (1-e)\begin{pmatrix}b\\ d\end{pmatrix}=0,\quad 
D\begin{pmatrix}b\\ d\end{pmatrix}=\bar\del e\begin{pmatrix}b\\ d\end{pmatrix}=-(1-e)\extd(1-e) \begin{pmatrix}b\\ d\end{pmatrix}\]
using the comutations above for $D$ and the relations (\ref{deld}) to find $\del e=e\extd e$ and $\bar\del e=-(1-e)\extd(1-e)=(\extd e)e$.

The map $\C_q[S^2]^{\oplus 2}\to \CS_+$ is given by $(f,g)\mapsto fb+gd=(f,g)e\begin{pmatrix}b\\ d\end{pmatrix}$. Simlarly for $\CS_-$ with $1-e$. Since given by complementary projectors, we see that together these maps trivialise $\CS$. Let us write this explicitly as the
combined map
\[ (f,g)\mapsto (f,g)\begin{pmatrix}a+\lambda b\\ c+\lambda d\end{pmatrix}\]
viewed in $\CS\subset \C_q[SL_2]$. The $a,c$ have different degrees from $b,d$ so any nonzero value of $\lambda$ will do and yield an isomorphism. Noting that
\[ \del e \begin{pmatrix}b\\ d\end{pmatrix}=\bar\del e \begin{pmatrix}a\\ c\end{pmatrix}=0\]
allows us to combine the expressions for $D$ according to
\[ D((f, g)\begin{pmatrix}a+\lambda b\\ c+\lambda d\end{pmatrix})=(\extd f,\extd g)\bar\tens \begin{pmatrix}a+\lambda b\\ c+\lambda d\end{pmatrix}-(f,g)(e\extd e-(\extd e)e)\begin{pmatrix}a+\lambda b\\ c+\lambda d\end{pmatrix}.\]
In view of the isomorphism, one can view this as an operator on $(f,g)$ which is as stated with $\extd$ and matrix multiplication acting from the right (this is inevitable
in our conventions).

For the Dirac operator we use the Leibniz formula (\ref{Dslleib}) and set $\lambda=q^{-\h}$ so that $a+\lambda b,c+\lambda d$ are eigenvectors of $\Dsl$. Then
\[ \Dsl((f, g)\begin{pmatrix}a+\lambda b\\ c+\lambda d\end{pmatrix})=\gamma_{(\extd f,\extd g)}(\begin{pmatrix}a+\lambda b\\ c+\lambda d\end{pmatrix})+\lambda^{-1}(f, g)\begin{pmatrix}a+\lambda b\\ c+\lambda d\end{pmatrix}.\]
By similar methods as in the proof of Proposition~\ref{Dsl}, we write $\extd f=f_i\extd b_i$ and take these coefficients also for expansions of $\del f,\bar\del f$. Then
\[ \gamma_{(\extd f,\extd g)}(\begin{pmatrix}a+\lambda b\\ c+\lambda d\end{pmatrix})=
(f_i,g_i)b_i \begin{pmatrix}q^2 b +\lambda q^{-1} a\\ q^2d+\lambda q^{-1}c\end{pmatrix}
+(q f_0-qg_-,f_+)\begin{pmatrix}b\\ d\end{pmatrix}+\lambda (g_-,-q^{-1}f_++g_0)\begin{pmatrix}a\\ c\end{pmatrix}.\]
By inserting an $e$ in front of $\begin{pmatrix}b\\ d\end{pmatrix}$ and $(1-e)$ in
front of $\begin{pmatrix}a\\ c\end{pmatrix}$, we can write this as a certain operation
on $(f,g)$ followed by dot product with $\begin{pmatrix}a+\lambda b\\ c+\lambda d\end{pmatrix}$, which gives the $\Dsl$ stated. Here $\overset{\leftarrow}\del_i$ denotes the right-acting operator
which assigns our chosen coefficients $f_i$ to $f$ and $g_i$ to $g$. The 
last two terms of $\Dsl$ (the displayed matrix terms) can be combined as
\[ \lambda\begin{pmatrix}-q\lvec\del_0 b_0-[2]_q \lvec\del_+ b_+& q^3\lvec\del_0 b_-+q\lvec\del_+(1+[2]_qb_0\\ \lvec\del_0 b_+ +
\lvec\del_-(1+[2]_q b_0)& -q^2 [2]_q \lvec\del_- b_- -q\lvec\del_0 b_0\end{pmatrix}.\]
This part alone is not well-defined (recall that the chosen coefficients of $\extd f,\extd g$ are not unique) but the entire $\Dsl$ is well-defined and independent of the choice, as it must be 
since $\Dsl$ exists geometrically on the bundle $\CS$. \eproof

It remains to make some remarks about the classical limits. Setting $q=1$ we have
\[ \Dsl_{q=1}=1+\begin{pmatrix}\lvec\del_- b_--\lvec\del_+ b_+& \lvec\del_0 b_-+\lvec\del_+(1+2b_0)\\ \lvec\del_0 b_++\lvec\del_-(1+2 b_0)& \lvec\del_+ b_+-\lvec\del_- b_-\end{pmatrix}=1-\imath\sigma\cdot\lvec\del\times\lvec x\]
when we shift to usual $x,y,z$ coordinates related (say) as $b_\pm=\pm(x\pm\imath y)$ and $b_0=z-\h$, corresponding to a sphere of radius $1/2$ embedded in $\R^3$. Changing
variables to $\lvec\del =(\lvec\del_x,\lvec\del_y,\lvec\del_z)$ (so that 
$\lvec\del_\pm=\pm{1\over 2}(\lvec\del_x\mp\imath\lvec\del_y)$ and $\lvec\del_0=\lvec\del_z$) gives $\Dsl$ as stated in terms of
the usual Pauli matrices and the vector $\lvec x=(x,y,z)$ acting by right 
multiplication. This expression makes sense on functions on $\R^3$ but 
vanishes on functions that depend only on the radius, and hence descends to functions on
$S^2$. Incidentally, as is well-known, $\Dsl^2=\Dsl+\square$ where
\[ \square =(\lvec\del \cdot \lvec x)^2+\lvec\del \cdot \lvec x- 
(\lvec\del\cdot\lvec\del)(\lvec x\cdot\lvec x)\]
(an easy computation) and one may check that $\square$ defined on the 
ambient $\R^3$ again descends to $S^2$ and is indeed the classical limit
of the Laplace operator in Proposition~\ref{laplace} under the same change of
classical coordinates and $\lvec x\cdot \lvec x={1/4}$. Obviously all operations here may be reformulated
acting from the left as more usual. 

\appendix
\section{Geometrical $q$-Borel-Weil-Bott construction}

The Borel-Weil-Bott construction in classical representation
theory constructs irreducible representations of a compact Lie
group $G$ as follows. From the bundle $G\to G/T$ (where $T$ is the
maximal torus) and an irreducible representation $V$ of $T$ (given
by a character) we construct the associated bundle $E=G\times_T V$
over $G/T$. Its space of sections $\CE$ still carries a
representation of $G$ acting on $G/T$ from the left and lifted to
an action on $\CE$ by a certain connection. This space is,
however, much too big to be an irreducible representation. From
the point of view of `geometric quantization' one must choose a
polarization. A natural approach here is to use the complex
structure on $G/T$ and take the holomorphic sections $\CE^{\rm
hol}$. These are a much smaller space and form an irreducible
representation of $G$. We refer to \cite{Zhang} for an excellent
account of the classical situation from this point of view and of
the quantum case from a representation theoretic (but not really
geometric) point of view. On the other hand, in the course of
understanding the geometry of the $q$-sphere, we have now obtained
all the ingredients for the quantum group geometrical version of
this construction. We outline this application in this appendix.
The full details and generalisation from $\C_q[SL_2]$ to other
quantum groups will be addressed elsewhere.

Indeed, this remark is about the monopole associated bundles and
their covariant derivatives as we have already used for charges
$\pm 2,\pm 1$. Now we consider general $\CE_{-n}=(\C_q[SL_2]\tens
V)^{\C[t,t^{-1}]}$ where $V=\C.v$ is the right comodule defined by
$\Delta_Rv=v\tens t^{-n}$. This is 1-dimensional so that
$\CE_n=\C_q[SL_2]_{n}.v$, i.e. isomorphic to the degree $n$
component. As such $\CE_{-n}$ also carries the left coaction of
$\C_q[SL_2]$ given by restricting the coproduct since this
respects the degree (the coaction on degree 0 was already used in
(\ref{coactS})). Let $D$ be the usual covariant derivative for the
monopole connection. We say that a section $\sigma\in \CE_{-n}$ is
{\em holomorphic} if \eqn{Ehol}{ D\sigma\in \Omega^{1,0}\bar\tens
\CE_{-n}.} In other words, if we write $D=\vardelta+\bar\vardelta$
for the $\Omega^{1,0}$ and $\Omega^{0,1}$ parts then we require
$\bar\vardelta\sigma=0$.

\begin{propos} The space of holomorphic sections $\CE_{-n}^{\rm hol}$ of the
charge $-n$ $q$-monopole bundle contains the standard
$n+1$-dimensional corepresentation of $\C_q[SL_2]$.
\end{propos}
\proof Here $n\ge 0$. The standard $n+1$-dimensional
corepresentation of $\C_q[SL_2]$ corresponds to the
$q$-deformation of the standard $n+1$-dimensional irreducible
representation of $SL_2$ and is given by \eqn{carep}{\{ c^s a^t\
|\ s+t=n,\ s,t\ge 0\}.} These are all in the degree $n$ component
of $\C_q[SL_2]$ and form a left corepresentation via the
restriction of the coproduct. We verify that as such, they are
homomorphic. Indeed, as in the computation for Theorem~\ref{cov},
we compute
\[ \extd c^s=[s;q^2]c^{s-1}(ce^0+q de^+),\quad \extd
a^t=[t;q^2]a^{t-1}(ae^0+qbe^+)\] as easily proven by induction.
Then
\[ \extd(c^s a^t)=[n;q^2]c^s a^t
e^0+q^{t}c^{s-1}a^{t-1}(q[s;q^2]+[n;q^2]bc)e^+\] by the Leibniz
rule. Hence using (\ref{omegamon}) we have
\[ D(c^s a^t)=\extd(c^s
a^t)-c^sa^t\omega(t^n)=q^{t}c^{s-1}a^{t-1}(q[s;q^2]+[n;q^2]bc)e^+.\]
The expressions are slightly simpler when $s=0$ or $t=0$. Next, we
move $c^{s-1}a^{t-1}$ to the far right and use (\ref{eplus}) to
see that $D(c^s a^t)\in \Omega^{1,0}\bar\tens \CE_{-n}$ by the
same argument as in the proof of Theorem~\ref{cov}. Hence these
elements are holomorphic as claimed. \eproof

Conversely, if $f\in \C_q[S^2]$ and $\sigma\in \CE^{\rm hol}_{-n}$
then $D(f \sigma)=\extd f\bar\tens \sigma+f D(\sigma)$ means that
$f\sigma$ cannot be holomorphic unless $\bdel f=0$ which, as
explained in Section~4 (at least for generic $q$), means $f$ a
multiple of 1. This reminds us that $\CE^{\rm hol}_{-n}$ is indeed
a complex vector space but not a $\C_q[S^2]$-module. As such we
have seen that it contains  ${\rm span}\{c^s a^t\}$, which are
linearly independent over $\C$ and give the usual
$n+1$-dimensional corepresentation. On the other hand, since the
dimension of $\CE^{\rm hol}_{-n}$ classically is $n+1$, this
should also be true for generic $q$ (since all our structures
deform with classical dimensions). In this case, by dimensions,
$\CE^{\rm hol}_{-n}={\rm span}\{c^s a^t\}$ i.e. not only contains
but coincides with the $q$-deformed $n+1$-dimensional
corepresentation of $\C_q[SL_2]$. This outlines a geometric proof
of the $q$-Borel-Weil-Bott construction.

\bigskip
\baselineskip 14pt

\end{document}